\documentclass[12pt]{article}

\usepackage{theorem}
\usepackage{amssymb}
\usepackage{amsmath}
\usepackage{amsfonts}
\usepackage{times}

\usepackage[paper=a4paper, top=2.5cm, bottom=2.5cm, left=2.5cm, right=2.5cm]{geometry}

\def\qed{\hfill $\vcenter{\hrule height .3mm
\hbox {\vrule width .3mm height 2.1mm \kern 2mm \vrule width .3mm
height 2.1mm} \hrule height .3mm}$ \bigskip}

\def \Sph{\mathbb{S}^{n-1}}
\def \RR {\mathbb R}

\def \EE {\mathbb E}

\def \PP {\mathbb P}
\def \eps {\varepsilon}
\def \vphi {\varphi}

\def \cE {\mathcal E}

\newtheorem{theorem}{Theorem}[section]
\newtheorem{lemma}[theorem]{Lemma}

\newtheorem{proposition}[theorem]{Proposition}
\newtheorem{corollary}[theorem]{Corollary}
 {\theorembodyfont{\rmfamily}}

\long\def\symbolfootnotetext[#1]#2{\begingroup
\def\thefootnote{\fnsymbol{footnote}}\footnotetext[#1]{#2}\endgroup}

\begin{document}

\title{Dimensionality and the stability of the Brunn-Minkowski inequality}
\author{Ronen Eldan\textsuperscript{*}
  ~ and ~ Bo`az Klartag\textsuperscript{*}}
\date{}
\maketitle

\symbolfootnotetext[1]{School of Mathematical Sciences, Tel-Aviv
University, Tel Aviv 69978, Israel. Supported in part
 by the Israel
Science Foundation and by a Marie Curie Reintegration Grant from the
Commission of the European Communities. Email:
\{roneneld,klartagb\}@tau.ac.il.}

\begin{abstract}
We prove stability estimates for the Brunn-Minkowski inequality for
convex sets. As opposed to previous stability results, our estimates
improve as the dimension grows. In particular, we obtain a
non-trivial conclusion for high dimensions already when
$$ Vol_n \left(
\frac{K + T}{2} \right)  \leq 5 \sqrt{ Vol_n(K) Vol_n(T) }. $$ Our
results are equivalent to a {\it thin shell} bound, which is one of
the central ingredients in the proof of the central limit theorem
for convex sets.
\end{abstract}

\section{Introduction}

The Brunn-Minkowski inequality states, in one of its normalizations,
that
\begin{equation}
 Vol_n \left( \frac{K + T}{2} \right) \geq \sqrt{Vol_n(K) Vol_n(T)}
 \label{eq_928}
\end{equation}
for any compact sets $K, T \subset \RR^n$, where $(K + T) / 2 = \{
(x + y) / 2 ; x \in K, y \in T \}$ is half of the Minkowski sum of
$K$ and $T$, and where $Vol_n$ stands for the Lebesgue measure in
$\RR^n$. Equality in (\ref{eq_928}) holds if and only if $K$ is a
translate of $T$ and both are convex, up to a set of measure zero.

\medskip The literature contains various stability estimates for
the Brunn-Minkowski inequality, which imply that when there is
almost-equality in (\ref{eq_928}), then $K$ and $T$ are almost-translates of each other. Such estimates appear in Diskant
\cite{diskant}, in Groemer \cite{groemer}, and in Figalli, Maggi and
Pratelli \cite{FMP1, FMP2}. We recommend Osserman \cite{osserman}
for a general survey on the stability of geometric inequalities.

\medskip All of the stability results that we found in the literature share a common
feature: Their estimates deteriorate quickly as the dimension
increases. For instance, suppose that $K, T \subset \RR^n$ are
convex sets  with
\begin{equation}
 Vol_n(K) = Vol_n(T) = 1 \quad \text{and} \quad \quad Vol_n \left(
\frac{K + T}{2} \right)  \leq 5. \label{eq_1211} \end{equation} The
present stability estimates do not seem to imply much about the
proximity of $K$ to a translate of $T$ under the assumption
(\ref{eq_1211}). Only if the constant ``$5$'' in (\ref{eq_1211})
is replaced by something like $1 + 1/n$ or so, then the results of
Figalli, Maggi and Pratelli  \cite{FMP2} can yield meaningful
information. The goal of this note is to raise the possibility that
the stability of the Brunn-Minkowski inequality actually {\it
improves} as the dimension increases. In particular, we would like
to deduce from (\ref{eq_1211}) that
\begin{equation}
 \left| \frac{\int_K p(x - b_K ) dx}{\int_T p(x - b_T) dx} \, - 1 \, \right| \ll 1 \label{eq_1217} \end{equation}
 for a family of non-negative functions
$p$, when the dimension $n$ is high. Here, $b_K$ and $b_T$ denote
the barycenters of $K$ and $T$ respectively. Furthermore, in some
non-trivial cases we may conclude (\ref{eq_1217})  even when the
constant ``$5$'' in (\ref{eq_1211}) is replaced by an expression
that grows with the dimension, such as $\log n$ or $n^{\alpha}$ for
a small universal constant $\alpha > 0$.

\medskip In this note we take  the first steps towards a
dimension-sensitive stability theory of the Brunn-Minkowski
inequality. First, let us focus  on the simplest case in which
$p(x)$ in (\ref{eq_1217}) is a quadratic polynomial. In fact, we are
interested mainly in expressions related to the quadratic form
\begin{equation}  q_K(x) = \frac{1}{Vol_n(K)} \int_K \langle x, y \rangle^2 d y \, - \, \left( \frac{1}{Vol_n(K)} \int_K \langle x, y \rangle d y \right)^2
\quad \quad (x \in \RR^n) \label{eq_1232} \end{equation} where
$\langle \cdot, \cdot \rangle$ is the standard scalar product in
$\RR^n$. The {\it inertia form} of the bounded, open set  $K \subset
\RR^n $ is defined as
\begin{equation}  p_K(x) = \sup \left \{ \langle x, y \rangle^2 \, ; \, q_K(y) \leq
1 \right \}. \label{eq_344} \end{equation} Note that $p_K$ is a
positive-definite quadratic form in $\RR^n$. We say that $K \subset
\RR^n$ is isotropic when the barycenter of $K$ lies at the origin
and $q_K(x) = |x|^2 = \langle x, x \rangle$ for all $x$. In this
case, also $p_K(x) = |x|^2$.  It is easy to see that any bounded,
open set $K \subset \RR^n$ has an affine image which is isotropic.

 \medskip A convex body in $\RR^n$ is a bounded, open convex
set. For a convex body $K \subset \RR^n$ we denote by $\mu_K$ the
uniform probability measure on $K$. Our first stability result is as
follows:

\begin{theorem} \label{thmsec2}
Let $K, T \subset \RR^n $ be convex bodies and let $R \geq 1$. Assume
that
$$ Vol_n \left( \frac{K+T}{2} \right) \leq R \sqrt{Vol_n(K) Vol_n(T)}. $$ Let $p(x) = p_K(x)$ be
the inertia form of $K$ defined in
(\ref{eq_1232}) and (\ref{eq_344}). Then,
\begin{equation}
 \left| \, \frac{ \int_T p(x - b_T) d \mu_T(x)}{\int_K p(x - b_K) d \mu_K(x)} \ - \
 1 \,
\right| \leq C \frac{R^{\alpha_2}}{n^{\alpha_1}}. \label{eq_348}
\end{equation}
 Here $C, \alpha_1, \alpha_2
> 0$ are universal constants, and $b_K, b_T$ are the barycenters of $K,T$ respectively.
\end{theorem}

See Theorem \ref{thmsec6} below for explicit bounds on the
universal constants $\alpha_1, \alpha_2$ from Theorem \ref{thmsec2}.
 Our interest in the inertia form $p_K$ stems from the
{\it central limit theorem for convex sets}, see \cite{EK, K1} for
background reading. As we shall explain in Proposition \ref{thin_shell}
below, Theorem \ref{thmsec2} implies the bound
\begin{equation} \sigma_n \leq C n^{1/2- \alpha_1} \label{eq_353}
\end{equation} where $\sigma_n$ is the {\it thin shell} parameter
from \cite{EK1}, $C > 0$ is a universal constant and $\alpha_1 > 0$
is the constant from Theorem \ref{thmsec2}. In fact, Theorem
\ref{thmsec6} and (\ref{eq_1238}) below show that the inequality
(\ref{eq_353}) is essentially an equivalence. Consequently, the
universal constant $\alpha_1$ from Theorem \ref{thmsec2} is
intimately connected with the thin shell parameter $\sigma_n$. The question
of whether $\sigma_n$ is bounded by a universal constant is currently
one of the central problems in high-dimensional convex geometry.

\medskip Next, we address the task of finding a larger class of functions $p$ for
which bounds such as  (\ref{eq_1217}) hold true. Suppose that
$\mu_1$ and $\mu_2$ are two Borel probability measures on $\RR^n$. A
Borel probability measure $\gamma$ on $\RR^n \times \RR^n$ is called
a {\it coupling} of $\mu_1$ and $\mu_2$ if $(P_1)_*(\gamma) = \mu_1$
and $(P_2)_*(\gamma) = \mu_2$ where $P_1(x,y) = x$ and $P_2(x,y) =
y$. Here, $(P_i)_*(\mu)$ denotes the push-forward of $\mu$ under the
map $P_i$ for $i=1,2$. For two Borel probability measures $\mu_1$
and $\mu_2$ on $\RR^n$ and for $1 \leq p < \infty$, we set
$$ W_p \left( \mu_1, \mu_2 \right) = \inf_{\gamma} \left( \int_{\RR^n \times \RR^n}  |x - y|^p d \gamma(x,y) \right)^{1/p} $$
where the infimum runs over all couplings $\gamma$ of $\mu_1$ and
$\mu_2$. This is precisely the $L^p$ Monge-Kantorovich-Wasserstein
transportation distance between $\mu_1$ and $\mu_2$. See, e.g.,
Villani's book \cite{villani} for more information about this
metric. Note that
 for any $1$-Lipschitz function $\vphi:
\RR^n \rightarrow \RR$,
$$
 \left| \int_{\RR^n} \vphi(x) d \mu_1(x) - \int_{\RR^n} \vphi(x) d \mu_2(x) \right|  \leq W_1(\mu_1,
 \mu_2) \leq W_2(\mu_1, \mu_2).
$$
In fact, the assumption that $\vphi$ is $1$-Lipschitz may typically be weakened. For
instance, when $\vphi$
 is convex or concave, it is well-known that
\begin{equation}
 \label{eq_1100}
    \left| \int_{\RR^n} \vphi  d \mu_1 -  \int_{\RR^n} \vphi d \mu_2 \right|   \leq W_2(\mu_1, \mu_2) \cdot \sqrt{ \max \left \{ \int_{\RR^n} |\nabla \vphi|^2 d \mu_1, \int_{\RR^n} |\nabla \vphi|^2 d \mu_2 \right \} }.
\end{equation}
\begin{theorem} \label{thmsec6_half}
Let $K, T \subset \RR^n $ be convex bodies whose barycenters lie at
the origin and let $R \geq 1$. Suppose that
$$ Vol_n \left( \frac{K+T}{2} \right) \leq R \sqrt{Vol_n(K) Vol_n(T)}. $$
Assume that $K$ is isotropic. Then,
\begin{equation} \label{thmsec6eq2}
\frac{W_2( \mu_K, \mu_T)}{\sqrt{n}} \leq C n^{-1/4} \sqrt{\sigma_n}
R^{5/2} \leq \tilde{C} \frac{R^{5/2}}{n^{\alpha}},
\end{equation}
where $\alpha, C, \tilde{C} > 0$ are universal constants.
\end{theorem}

Theorem \ref{thmsec6_half} combined with the inequality
(\ref{eq_1100}) entails the bound (\ref{eq_1217}) in the case where,
for instance, $p(x) = \| x \|^q$ for various norms $\| \cdot \|$ in $\RR^n$,
$q \geq 0$ and $R \ll n^c$. Additionally, the estimate
(\ref{thmsec6eq2}) implies the non-trivial bound (\ref{eq_348}) via
(\ref{eq_1100}). We do not know  the optimal value of the
exponent $\alpha$ in Theorem \ref{thmsec6_half}. We know more in the
particular case of {\it unconditional convex bodies}.
 A convex
body in $\RR^n$ is said to be {\it unconditional} if $$
(x_1,\ldots,x_n) \in K \quad \Longleftrightarrow \quad (\pm x_1,\ldots,\pm
x_n) \in K $$ for all $(x_1,\ldots,x_n) \in \RR^n$ and for all possible
choices of signs. In other words, $K$ is invariant under coordinate
reflections. For unconditional convex bodies, Theorem
\ref{thmsec6_half} may be sharpened as follows:

\begin{theorem} \label{uncond}
Let $K, T \subset \RR^n$ be unconditional convex bodies, and let $R \geq
1$. Assume that $K$ is isotropic and that
$$ Vol_n \left( \frac{K+T}{2} \right) \leq R \sqrt{Vol_n(K) Vol_n(T)}. $$
 Then
\begin{equation}  W_2( \mu_K, \mu_T) \leq C (R-1)^{5/2} \log n, \label{eq_305}
\end{equation}
 where $C
> 0$ is a universal constant. \label{thm_1455}
\end{theorem}

Thus, in the unconditional case, the exponent $\alpha$ from Theorem
\ref{thmsec6_half} is essentially $1/2$, up to logarithmic factors.
When substituting $\vphi(x) = |x|^2$ in (\ref{eq_1100}) and using
(\ref{eq_305}), we conclude that for any $K, T \subset \RR^n$ as in
Theorem \ref{uncond},
\begin{equation}  \left| \int_K |x|^2 d \mu_K - \int_T |x|^2 d \mu_T \right| \leq C \sqrt{n} \cdot \log n  \cdot (R - 1)^5 \label{eq_1125}
\end{equation}
(in order to use (\ref{eq_1100}) we also need a crude estimate for
$\int_T |x|^2 d \mu_T(x)$, hence we applied Corollary \ref{boundedvar} to
obtain such an estimate). In view of (\ref{eq_1125}) and Proposition
\ref{thin_shell} below, we match (up to logarithmic factors) the
best bounds for the width of the thin spherical shell for
unconditional convex bodies proven in \cite{uncond}.

\medskip
The structure of the remainder of this note is as follows: In the
next section we establish some well-known facts about
one-dimensional log-concave measures. In Section 3 we prove Theorem \ref{thmsec2} and in Section 4 we prove Theorem \ref{thmsec6_half}.
Section 5 is dedicated to attaining some inequalities related to one-dimensional transportation of measure. In Section 6, using these inequalities, we prove Theorem \ref{uncond}.

\medskip Throughout this note, we write $c, C, \tilde{c}$ etc. for
various positive universal constants, whose value may change from
one line to the next. We usually use upper-case $C$ to denote
universal constants that we consider  ``sufficiently large'', and
lower-case $c$ to denote universal constants that are ``sufficiently
small''. We write $\log$ for the natural logarithm. By ``measurable'' we always mean Borel-measurable.

\section{Background about log-concave densities on the line}

In this section we recall some facts, all of which are well-known to
experts, about log-concave densities. A function $\rho: \RR^n
\rightarrow [0, \infty)$ is log-concave if for any $x, y \in \RR^n$,
$$ \rho \left( \lambda x + (1 - \lambda) y \right) \geq
\rho(x)^{\lambda} \rho(y)^{1-\lambda} \ \ \ \ \ \text{for all} \ 0 <
\lambda < 1. $$
 A probability measure or a random variable are called log-concave
if they posses a log-concave density. Let $\mu$ be a log-concave
probability measure on $\RR$, whose log-concave density is denoted
by $\rho: \RR \rightarrow [0, \infty)$. Write
$$ \Phi(t) = \mu \left( (-\infty, t] \right) = \int_{-\infty}^t \rho(s) ds \ \ \ \ \ \ \ \ \ \ (t \in \RR). $$
A nice characterization of log-concavity that we learned from Bobkov
\cite{B2} is that $\mu$ is log-concave if and only if the function
$$ t \mapsto \rho(\Phi^{-1}(t)) \ \ \ \ \ \ \ \ \ \ \ \ \ t \in
[0,1] $$ is a concave function. This characterization lies at the
heart of the proof of the following Poincar\'e-type inequality which
appears as Corollary 4.3 in Bobkov \cite{B1}:

\begin{lemma}  \label{spectral_gap}
Let $\mu$ be a log-concave probability measure on the real line, and
set
$$ Var(\mu) = \int x^2 d\mu(x) - \left( \int x d \mu(x)
\right)^2 $$ for the variance of $\mu$. Then for any smooth function
$f$ with $\int f d \mu = 0$,
$$
\int_{\RR} f^2(t) d \mu(t) \leq 12 Var(\mu) \int_{\RR} |f'(t)|^2 d
\mu(t).
$$
\end{lemma}

Further information about log-concave densities on the line is
provided by the following standard lemma.

\begin{lemma} Let $f: \RR \rightarrow [0, \infty)$ be a log-concave
probability density. Denote $b = \int x f(x) dx$, the barycenter of
the density $f$, and let $\sigma^2$ be the variance of the
random variable whose density is $f$. Then, for any $t \in \RR$,
\begin{enumerate}
\item[(a)] $\displaystyle  f(t) \leq \frac{C}{\sigma} \exp(-c |t - b| / \sigma)$;
and
\item[(b)] If $|t - b| \leq c \sigma$, then $ \displaystyle f(t) \geq
\frac{c}{\sigma}$.
\end{enumerate}
Here, $c, C > 0$ are universal constants. \label{lem_310}
\end{lemma}

\textbf{Proof:} Part (a) is the content of Lemma 3.2 in Bobkov
\cite{B3}. In order to prove (b), we show that for some $t_0 \geq b
+ c_0 \sigma$,
\begin{equation}
f(t_0) \geq 1 / (10 C_1 \sigma) \label{eq_245}
\end{equation}
with $c_0 = 1 / (10 C)$, $C_1 =  c^{-1} \log(10 C / c)$ where here
$c, C$ are the constants from part (a). Indeed, if there is no such
$t_0$, then from (a),
$$ \int_b^{\infty} f(t) dt \leq \int_b^{b+c_0 \sigma} \frac{C}{\sigma} dt + \int_{b+c_0 \sigma}^{b+C_1 \sigma } \frac{dt}{10 C_1 \sigma} + \int_{b +
C_1 \sigma}^{\infty} \frac{C}{\sigma} \exp(-c |t - b| / \sigma) dt
\leq \frac{3}{10} < \frac{1}{e},
$$
in contradiction to Gr\"unbaum's inequality (see, e.g., \cite[Lemma
3.3]{B3}). By symmetry, there exists some $t_1 \leq b - c_0 \sigma$
with
$$ f(t_1) \geq  1 / (10 C_1 \sigma). $$
From log-concavity, $f(t) \geq 1 / (10 C_1 \sigma)$ for $t \in [t_1,
t_0]$, and (b) is proven since $[t_1, t_0] \supseteq [b - c_0 \sigma, b + c_0 \sigma]$. \qed

\medskip  The
following lemma is essentially a one-dimensional, functional version of Theorem \ref{thmsec2}.   The Lemma states, roughly, that if the supremum-convolution of
two log-concave probability densities has a bounded integral,
then their respective variances cannot be too far from each other.

\begin{lemma} \label{onedimvar}
Let $X,Y$ be random variables with corresponding densities $f_X,
f_Y$ and variances   $\sigma_X^2, \sigma_Y^2$. Assume that $f_X$ and
$f_Y$ are log-concave. Define
\begin{equation} \label{supconv}
h(t) = \sup_{s \in \RR} \sqrt{f_X(t+s) f_Y(t-s)},
\end{equation}
a supremum-convolution of $f_X$ and $f_Y$. Then,
$$
\int_{\RR} h(t) dt \geq c \sqrt {\max \left \{
\frac{\sigma_X}{\sigma_Y}, \frac{\sigma_Y}{\sigma_X} \right \}}
$$
where $c > 0$ is a universal constant.
\end{lemma}

\textbf{Proof:} The function $h$ is clearly measurable (it is even log-concave). It follows from Lemma \ref{lem_310}(b) that
there exist intervals $I_X, I_Y$ such that
$$
Length(I_X) \geq c \sigma_X, ~~ Length(I_Y) \geq c \sigma_Y
$$
and,
$$
f_X(t) \geq \frac{c}{\sigma_X}, ~~ \forall t \in I_X ~~;~~ f_Y(s)
\geq \frac{c}{\sigma_Y}, ~~ \forall s \in I_Y.
$$
Combining this with (\ref{supconv}), we learn that there exists an
interval $I_Z$ with $Length(I_Z) \geq c (\sigma_X + \sigma_Y) / 2$
such that,
$$
h(t) \geq \frac{c}{\sqrt {\sigma_X \sigma_Y} }, ~~ \forall t \in
I_Z.
$$
This implies,
$$
\int_{\RR} h(t) dt \geq \int_{I_Z} h(t) dt \geq \frac{c^2}{2}
\frac{\sigma_X + \sigma_Y}{\sqrt{\sigma_X \sigma_Y}} \geq
\frac{c^2}{2} \sqrt {\max \left \{ \frac{\sigma_X}{\sigma_Y},
\frac{\sigma_X}{\sigma_Y} \right \}}
$$
which completes the proof. \qed

\medskip Recall the definition (\ref{eq_1232}) of the inertia form
$q_K(x)$ associated with a convex body $K \subset \RR^n$.  As a
corollary of Lemma \ref{onedimvar}, we have

\begin{corollary} \label{boundedvar}
Let $R > 1$ and let $K,T \subset \RR^n$ be convex bodies such that
$$ Vol_n \left(\frac{K+T}{2} \right) \leq R \sqrt{Vol_n(K) Vol_n(T)}. $$ Then,
\begin{equation}
\frac{1}{C R^4} q_K(x) \leq q_T(x) \leq C R^4 q_K(x) \quad \quad
\quad \text{for all} \ x \in \RR^n \label{eq_348_}
\end{equation} where $C
> 0$ is a universal constant.
\end{corollary}
\textbf{Proof:}  Fix a unit vector $\theta \in \RR^n$.  Let $\tilde
X, \tilde Y$ be random vectors distributed uniformly  on $K, T$
respectively, and define $X = \langle \tilde{X}, \theta \rangle$ and
$Y = \langle \tilde{Y}, \theta \rangle$. Observe that
$$ q_K(\theta) = Var(X), \quad q_T(\theta) = Var(Y). $$
In order to prove (\ref{eq_348_}), it suffices to show that
\begin{equation}
\max \left \{ \frac{Var(X)}{Var(Y)}, \frac{Var(Y)}{Var(X)} \right \}
\leq C R^4. \label{eq_350}
\end{equation}
Denote the respective densities of $X, Y$ by $f_X, f_Y$. The
Pr\'ekopa-Leindler theorem (see, e.g., the first pages of Pisier
\cite{pisier}) implies that $f_X$ and $f_Y$ are log-concave.
Furthermore, using the Pr\'ekopa-Leindler theorem again we derive,
\begin{equation} \label{supconv2}
Vol_n \left ( \frac{K+T}{2} \right ) \geq \int_\RR \sup_{s \in
\mathbb{R}}  \sqrt{f_X(t-s) Vol_n(K) f_Y(t+s) Vol_n(T)} dt.
\end{equation}
Hence,
$$
\int_{\RR} \sup_{s \in \mathbb{R}}  \sqrt{f_X(t-s) f_Y(t+s)} dt \leq
R.
$$
Plugging this into lemma \ref{onedimvar} we deduce (\ref{eq_350}).
\qed

\medskip \textbf{Remark.} Let $K, T, R$ be as in Corollary \ref{boundedvar}
and let $\tilde
X, \tilde Y$ be the random vectors  distributed uniformly on $K, T$ respectively.
Corollary \ref{boundedvar} states that
\begin{equation}
\frac{1}{C R^4} Cov(\tilde{X}) \leq Cov(\tilde{Y}) \leq C R^4 Cov(\tilde{X})
\label{eq_1619}
\end{equation}
in the sense of symmetric matrices, where $Cov(\tilde{X})$ is the covariance matrix of $\tilde{X}$.
Furthermore, we do not have to assume that $\tilde{X}, \tilde{Y}$ are distributed uniformly
in a convex body. The estimate (\ref{eq_1619}) holds true whenever $\tilde{X}, \tilde{Y}$ have log-concave densities
$f_{\tilde{X}}, f_{\tilde{Y}}$ with
$$ R = \int_{\RR^n} \left( \sup_{y \in \RR^n} \sqrt{f_{\tilde{X}}(x+y) f_{\tilde{Y}}(x-y)} \right) dx. $$

\medskip
Next, for a measure $\mu$ and measurable sets $A,B$ with $0 <
\mu(A) < \infty$ define
$$
\mu|_A(B) = \frac{\mu(A \cap B)}{\mu(A)}.
$$
Thus the probability measure $\mu|_A$ is the conditioning of $\mu$ to the set $A$. Clearly, for a
log-concave measure $\mu$ and an interval $I$, the measure $\mu|_I$
remains log-concave.
\begin{lemma} Let $\mu$ be a log-concave probability measure on $\RR$.
Then for any two intervals $J_1 \subseteq J_2 \subset \RR$,
$$ Var(\mu|_{J_1}) \leq Var(\mu|_{J_2}) $$
(the ``intervals'' may also include rays, or the entire line: Any
convex set in $\RR$). \label{lem_1213}
\end{lemma}

\textbf{Proof:} It is enough to prove the lemma for $J_1, J_2$ being
rays. Denote by $I$ the interior of the support of $\mu$, and by
$\rho$ the density of $\mu$. Abbreviate $\Phi(t) = \mu \left (
(-\infty, t ] \right ), \ \mu_t = \mu|_{(-\infty, t]}$ and set
$$ e(t) = \int_{\RR} x d \mu_t(x), \ \ \ \ \ \ \ \ v(t) = Var(\mu_t) = \int_{\RR} x^2 d
\mu_t(x) - e^2(t) \ \ \ \ \ \ \ \ \ \ \ \ \ (t \in I). $$ Then for any
 $t \in I$,
$$ e^{\prime}(t) = \frac{\rho(t)}{\Phi(t)} \left(  t - e(t) \right), \ \ \ \ \ \ \
v^{\prime}(t) =  \frac{\rho(t)}{\Phi(t)} \left( (t - e(t))^2 - v(t)
\right). $$ To prove the lemma, it suffices to show that
$v^{\prime}(t) \geq 0$ for any $t$, or equivalently, that
$$ Var(\mu_t) - (t - \EE \mu_t)^2 = v(t) - (t - e(t))^2 \leq 0 \ \ \ \ \ \text{for all} \ t \in I. $$
This is equivalent to demonstrating  that for any log-concave random
variable $X$ such that $X \geq 0$ almost surely,
one has $Var[X] \leq \left( \EE [X] \right)^2$. This follows immediately from Borell
\cite[Lemma 4.1]{borell}, see also Lov\'asz and Vempala \cite[Lemma
5.3(c)]{LV}. \qed

\section{Deriving a stability estimate from the central limit theorem  for convex sets }

In this section we prove Theorem \ref{thmsec2}.
The main ingredient we use is the central limit theorem for convex sets, proven initially in
\cite{K1}. It states that for any isotropic convex body $K \subset
\RR^n$, and for ``most'' subspaces of a small enough dimension, the
marginal of $\mu_K$ is approximately Gaussian. Below we use a
pointwise version of this theorem, proven in \cite{EK}, which shows
that there exists a subspace of dimension $n^{c}$, where $c>0$ is some universal constant, on which the
marginals of both $K$ and $T$ are both approximately Gaussian
density-wise. The Pr\'ekopa-Leindler inequality then implies that
the marginal of $(K+T)/2$ on the same subspace is pointwise greater
than the supremum-convolution of the respective marginals of $K$ and
$T$. Therefore, the density of the marginal of $(K+T)/2$  must be greater than the supremum-convolution of two
densities which are both approximately Gaussian, but typically have
different covariances.

\medskip
A second ingredient will be a calculation
which shows that the integral of the supremum-convolution of two
Gaussian densities whose covariance matrix is a multiple of the
identity, becomes very large when their respective covariances are not
close to one another. This will imply that when $Vol_n((K+T)/2)$ is
not large, the covariance matrices of both marginals are roughly the
same multiple of the identity. Therefore the inertia forms of $K$
and $T$ must have had roughly the same trace (the trace of the
matrix will determine the multiple of the identity).

\medskip We write $G_{n, \ell}$ for the Grassmannian of all
$\ell$-dimensional subspaces in $\RR^n$, and $\sigma_{n, \ell}$
stands for the Haar probability measure on $G_{n, \ell}$. A random
vector $X$ in $\RR^n$ is centered if $\EE X = 0$ and is isotropic if
its covariance matrix is the identity matrix. For a subspace $E
\subseteq \RR^n$ we write $Proj_E$ for the orthogonal projection
operator onto $E$ in $\RR^n$. Furthermore, define $\gamma_{k,\alpha} (x) = (2 \pi \alpha^2)^{-k/2} \exp(
-\frac{|x|^2}{2 \alpha^2} )$ the centered Gaussian density in $\RR^k$
with covariance $\alpha^2$, and abbreviate $\gamma_k(x) =
\gamma_{k,1}(x)$. The main result of \cite{EK} reads as follows:

\begin{theorem} \label{pointwise}
Let $X$ be a centered, isotropic random vector in $\RR^n$ with a
log-concave density. Let $1 \leq \ell \leq n^{c_1}$ be an integer.
Then there exists a subset $\cE \subseteq G_{n,\ell}$ with
$\sigma_{n,\ell}(\cE) \geq 1 - C \exp( -n^{c_2} )$ such that for any
$E \in \cE$, the following holds: Denote  by $f_E$ the log-concave
density of the random vector $Proj_E(X)$. Then,
\begin{equation}
\left | \frac {f_E(x)}{\gamma_\ell(x)} - 1 \right | \leq
\frac{C}{n^{c_3}} \label{eq_1241}
\end{equation}
for all $x \in E$ with $|x| \leq n^{c_4}$. Here, $C, c_1,c_2,c_3,c_4 > 0$ are universal
constants.
\end{theorem}

It can be seen directly from the proof in \cite{EK} that the
constants in Theorem \ref{pointwise} may be selected to be
$c_1,c_2,c_3=\frac{1}{30}, c_4=\frac{1}{60}, C=500$. Other
 constants would imply different universal constants in Theorem
 \ref{thmsec2}. We shall need the following elementary lemma:

\begin{lemma} For any $a > 0$, $$  \frac{1+a}{2 \sqrt{a}} \geq 1 + c \cdot \min \{ (\alpha-1)^2,  1 \}, $$
for $\alpha = \sqrt{1/a}$ and also for $\alpha = a$, where $c > 0$ is a universal constant.
 \label{lem_1457}
\end{lemma}

\textbf{Proof:} First we prove the lemma for $\alpha = a$. Note that for $0 < a \leq 4$,
$$ \frac{1+a}{2 \sqrt{a}} = 1 + \frac{1 - 2\sqrt{a} + a}{2 \sqrt{a}}
= 1 + \frac{(\sqrt{a} - 1)^2}{2 \sqrt{a}} = 1 + \frac{(a-1)^2}{2
\sqrt{a} (\sqrt{a}+1)} \geq 1 + \frac{(a-1)^2}{12}, $$
while for $a > 4$ we may write
$$ \frac{1+a}{2 \sqrt{a}} = 1 + \frac{(\sqrt{a} - 1)^2}{2 \sqrt{a}} \geq  1 + \frac{\sqrt{a} - 1}{2 \sqrt{a}} \geq 1 + \frac{\sqrt{a} / 2}{2 \sqrt{a}} = 1 + \frac{1}{4}.  $$
The case where $\alpha = \sqrt{1/a}$ follows as $\min \{ (\sqrt{1/a}-1)^2,  1 \} \leq 10 \min \{ (a-1)^2,  1 \}$.
 \qed

The following lemma is the second ingredient in our proof of Theorem \ref{thmsec2} described above.
The essence of the lemma is that the integral of the supremum-convolution of two spherically-symmetric Gaussian densities
must be quite large when the covariances are not close to each other.

\begin{lemma} \label{gaussiansupconv}
Let $k \in \mathbb{N}~$ and $A, B, \alpha >0$. Let $f,g,h: \RR^k \to
[0, \infty)$ satisfy
$$
h(x) \geq \sup_{y \in \mathbb{R}^k}  \sqrt{f(x-y) g(x+y)}, ~~~~ \forall x \in \RR^k
$$
and suppose that,
$$
f(x) \geq A \gamma_{k,1}(x)
$$
whenever $|x| \leq 10 \sqrt k$, and that
$$
 g(x) \geq B \gamma_{k,\alpha}(x), ~~
$$
whenever $|x| \leq 10 \alpha \sqrt k $. Assume that $h$ is measurable. Then,
\begin{equation}
\int_{\RR^k} h(x) dx \geq \frac{1}{2} \sqrt{A B } \left( 1 + c \cdot \min \{ (\alpha-1)^2,  1 \} \right)^{k/4}, \label{eq_1709}
\end{equation}
where $c > 0$ is a universal constant.
\end{lemma}

\textbf{Proof:} By homogeneity, we may  assume that $A = B = 1$.
Denote $a = 1 / \alpha^{2}$. Fix a unit vector $\theta \in \RR^n$
and $t > 0$. Then for any $s \in \RR$ with $ |s + t| \leq 10
\sqrt{k}$ and $|s - t| \leq 10 \alpha \sqrt{k}$,
\begin{equation}
h(t \theta) \geq \sqrt {f((t+s) \theta) g ((t-s) \theta)} \geq \left
(\frac{\sqrt a}{2 \pi} \right )^{k/2} \exp \left(-\frac{1}{4}
((t+s)^2 + a(t-s)^2) \right). \label{eq_1346}
\end{equation}
We would like to find $s$ which maximizes the right-hand side in
(\ref{eq_1346}). We select $s = t (a - 1) / (a + 1)$ and verify that
when $|t| < 5 \sqrt{(1 + a)k / a}$ we have $|s + t| \leq 10
\sqrt{k}$ and $|s - t| \leq 10 \alpha \sqrt{k}$. We conclude that
for any $|t| < 5 \sqrt{(1 + a)k / a}$,
$$ h(t \theta) \geq \left
(\frac{\sqrt a}{2 \pi} \right )^{k/2} \exp \left(- t^2 a / (1 + a)
\right). $$ Consequently,
\begin{align*}
 \int_{\RR^k} h(x) dx & \geq \left
(\frac{\sqrt a}{2 \pi} \right )^{k/2} \int_{5 \sqrt{(1 + a)k / a}
B_2^k} \exp \left(- \frac{a |x|^2}{1 + a} \right) dx \\ & = \left
(\frac{1 + a }{4 \pi \sqrt a} \right )^{k/2} \int_{ \sqrt{50 k}
B_2^k} \exp \left(- \frac{|x|^2}{2} \right) dx \geq \frac{1}{2}
\left (\frac{1 + a }{2 \sqrt a} \right )^{k/2},
\end{align*} where $B_2^k = \{ x \in \RR^k ; |x| \leq 1 \}$, and
where we utilized the fact that $$ \PP( |Z|^2 \geq 50 k ) \leq \EE |Z|^2
/ (50 k) = \frac{1}{50} <  1/2 $$ when $Z$ is a standard Gaussian in
$\RR^k$. All that remains is to apply Lemma \ref{lem_1457}. \qed

The following lemma combines Theorem \ref{pointwise} with the estimate
we have just proved. For a probability density $g$ on $\RR^n$ we write
$Cov(g)$ for the covariance matrix of the random
vector with density $g$. We similarly define $Cov(\mu)$ for a probability measure $\mu$
on $\RR^n$.

\begin{lemma} \label{lemsec5}
Let $f,g$ be log-concave probability densities on $\RR^n$ such that $f$ is isotropic. Let $\{\lambda_i\}_{i=1}^n$ be the eigenvalues of $Cov(g)$, repeated according to their multiplicity. Denote
$$
R = \int_{\RR^n} \sup_{y \in \RR^n} \sqrt{ f(x+y) g(x-y)} dx.
$$
Then, for $0 < \delta < 1$,
$$
\# \{ i ~;~ |\lambda_i - 1| \geq \delta \} \leq C \left( \frac{\log (2R)}{\delta} \right)^{C_1}
$$
for some universal constants $C, C_1 > 1$.
\end{lemma}

\textbf{Proof:} Clearly, we may  assume that the sequence $\{ \lambda_i \}$ is non-decreasing.
Translating $g$, we may assume that the barycenter of $g$ is
at the origin.
Let $X$ and $Y$ be random vectors that are distributed according to the laws $f,g$, respectively.
Fix $0 < \delta < 1$. Consider the subspace $E$
spanned by $\{ e_i ; \lambda_i - 1 \geq \delta \}$, where $\{e_i\}$ is an orthonormal basis of eigenvectors corresponding
the the eigenvalues $\{ \lambda_i \}$.
Denote $d = \dim E$ and assume that $d \geq 2$.
Since the $\lambda_i$'s are in increasing order, the subspace $E$ has the form,
$$
E = span \{e_i, i \geq i_0 \}
$$
for some $1 \leq i_0 \leq n$. Write $j_0 = \left \lfloor \frac{n - i_0}{2} \right \rfloor$ and
$V^2 = \lambda_{i_0 + j_0}$.
Now, fix $1 \leq j \leq j_0$. Define,
$$
v_j(\theta) = \theta e_{i_0 + j_0 + j} + \sqrt{1 - \theta^2} e_{i_0 + j_0 - j}.
$$
Inspect the function $f(\theta) = \langle Cov(g) v_j(\theta), v_j(\theta) \rangle$. We have $f(0)
= \lambda_{i_0 + j_0 - j} \leq V^2$ and $f(1) = \lambda_{i_0 + j_0 + j} \geq V^2$. By continuity, there
exists a certain $0 \leq \theta_j \leq 1$ for which
\begin{equation} \label{vjworks}
\langle Cov(g) v_j(\theta_j), v_j(\theta_j) \rangle = V^2.
\end{equation}
Denote
$$
F = span \left \{v_j(\theta_j)~| ~~ 1 \leq j \leq j_0 \right \}.
$$
Equation (\ref{vjworks}) and the fact that $e_1,\ldots,e_n$ are orthonormal eigenvectors imply that for every
$v \in F$, one has $\langle Cov(g) v, v \rangle = V^2$.
Moreover, $\dim F = j_0 \geq \frac{1}{2} d - 1$.
We now apply Theorem \ref{pointwise} which claims that if $d \geq C$,
then there exists a subspace $G \subset F$ with $\dim G = \lfloor
d^{1/40} \rfloor$ such that
$$
\tilde f (x) \geq \frac{1}{2} \gamma_{k, 1} (x), ~~
\tilde g (y) \geq \frac{1}{2} \gamma_{k, V} (y)
$$
for all $x$ with $|x| \leq 10 d^{1/80}$ and for all $|y| \leq 10 V d^{1/80}$, where $\tilde f$ and $\tilde g$ are the densities of $Proj_G(X), Proj_G(Y)$ respectively. Next, we use Lemma \ref{gaussiansupconv} to attain
$$
\int_{G} \sup_{y \in G}  \sqrt{\tilde f (x-y) \tilde g(x+y)}
dx  \geq \frac{1}{4} (1 + c \cdot \min \{ (V-1)^2, 1 \})^{\dim G / 4}.
$$
On the other hand, we may use the Prekop\'a-Leindler inequality as in (\ref{supconv2})
above, and deduce that
$$
\int_{G} \sup_{y \in G}  \sqrt{\tilde f(x-y) \tilde g (x+y)}
dx \leq R.
$$
Consequently, under the assumption that $d \geq C$,
\begin{equation}  \min \left \{ (V - 1)^2, 1 \right \} \leq C \log (2 R) / \dim(G).
\label{eq_1549}
\end{equation}
Since $V \geq \sqrt{1 + \delta} \geq 1 + \delta/3$, we conclude
$$
\# \{ i ~;~ \lambda_i - 1 \geq \delta \} \leq C \left( \frac{\log (2R)}{\delta} \right)^{C_1}.
$$
By repeating the
argument, with the subspace $\{ e_i ; \lambda_i - 1 \leq -\delta \}$
replacing the subspace $E$, we conclude the proof. \qed

\medskip
\textbf{Proof of Theorem \ref{thmsec2}:}
By applying affine transformations to both $K$ and $T$, we can assume that both bodies
have the origin as their barycenter, and that
 $p_K(x) =
|x|^2$ while $p_T(x) = \sum_i x_i^2 / \lambda_i$. By Lemma \ref{lemsec5},
\begin{equation}  \# \left \{ i ; |\lambda_i - 1| \geq \delta \right \} \leq C \left( \frac{\log (2R)}{\delta} \right)^{C_1},
\label{eq_1127} \end{equation}
for any $0 < \delta < 1$. Since $\lambda_i \leq C R^4$ for all $i$, as follows
from Corollary \ref{boundedvar}, then
\begin{equation} \frac{1}{n} \sum_{i=1}^n (\lambda_i - 1)^2 \leq  \frac{C}{n} \int_0^1  \min \left\{ n, \left( \frac{\log (2R)}{\delta} \right)^{C_1} \right \} d \delta
+ \frac{\tilde{C} ( \log (2 R))^{C_1} R^4}{n}  \leq C \frac{R^{\alpha_2}}{n^{\alpha_1   }} \label{eq_1136}
\end{equation} where $C, \alpha_1, \alpha_2 > 0$ are universal constants.  To obtain (\ref{eq_348}), note that
\begin{equation}
 \left| \, \frac{ \int_T p_K(x - b_T) d \mu_T(x)}{\int_K p_K(x - b_K) d \mu_K(x)} \ - \
 1 \,
\right| = \frac{1}{n} \left| \sum_{i=1}^n (\lambda_i - 1) \right| \leq \sqrt{  \frac{1}{n} \sum_{i=1}^n (\lambda_i - 1)^2 }. \label{eq_1234}
\end{equation}
\qed

\medskip \textbf{Remark:} When $K$ in Theorem \ref{thmsec2} is isotropic, we actually prove in (\ref{eq_1136})
 that \begin{equation} \| Cov(\mu_K) - Cov(\mu_T) \|_{HS}^2 \leq C R^{\alpha_2} n^{1-\alpha_1}, \end{equation} where $\| A \|_{HS}^2 = Trace(A^t A)$ is the square of the Hilbert-Schmidt norm of the matrix $A$.

\section{Obtaining stability estimates using a transportation argument}

The goal of this section is to prove Theorem \ref{thmsec6_half} and to
obtain some quantitative estimates for the exponents
from Theorem \ref{thmsec2}.
We begin with several core definitions which will be used in the proof. For two measurable functions $f,g: \RR^n \to [0, \infty)$, denote by $H_\lambda (f,g)$ the supremum-convolution of the two functions, hence,
\begin{equation}
\label{eq_1145}
H_\lambda (f,g)(x) := \sup_{y \in \RR^n} f^{1 - \lambda}(x + \lambda y) g^{\lambda} (x - (1 - \lambda) y).
\end{equation}
The function
$$ (\lambda, x) \mapsto H_{\lambda}(f,g)(x) $$
is log-concave in $[0,1] \times \RR^n$. We define
$$
K_\lambda (f,g) = \int_{\RR^n} H_\lambda (f,g)(x) dx
$$
the integral over  a subspace, and
$$
K(f,g) = \int_0^1 K_\lambda (f,g) d \lambda,
$$
the entire integral. Next, we write
$$
b(f,g) = \frac{1}{K(f,g)} \int_{\RR^n} \int_0^1 x H_\lambda(f,g)(x) d \lambda dx,
$$
the barycenter of $\int_0^1 H_{\lambda}(f,g)(x) d \lambda$. For $x \in \RR^n$ we write
$x \otimes x = (x_i x_j)_{i,j=1,\ldots,n}$, an $n \times n$ matrix. Set
\begin{equation} \label{defmatd}
D(f,g) = \frac{1}{K(f,g)}\int_{\RR^n} \int_0^1 \left( x \otimes x \right) H_\lambda(f,g)(x + b(f,g)) d \lambda dx,
\end{equation}
the covariance matrix. Finally, we normalize this density by defining
$$
L(f,g)(\lambda, x) = \frac{1}{K(f,g)} \sqrt{\det D(f,g)} \cdot H_\lambda(f,g)(D^{1/2} x + b(f,g))
$$
and
$$
l(f,g)(x) = \int_0^1 L(f,g) (\lambda, x) d \lambda,
$$
the marginal of $L(f,g)$ with respect to the axis $\lambda$. Note that by the Pr\'ekopa-Leindler inequality, $l(f,g)$ is an isotropic log-concave probability density in $\RR^n$.

\medskip The results of this section rely on the so-called \emph{Brenier map} between two given log-concave measures. Given two smooth log-concave probability densities $f,g$ on $\RR^n$, one may consider the Monge-Amp\`{e}re equation,
$$
\det(Hess \varphi) = \frac{g \circ \nabla \varphi}{f}.
$$
A theorem of Brenier asserts that a convex solution to the above
equation on the domain $Supp(f) = \{ x ; f(x) > 0 \}$  exists. The
regularity theory developed by Caffarelli implies that the convex
function $\varphi$ is smooth. For precise definitions and
properties, see \cite{villani}. The map $F = \nabla \varphi$ pushes
forward the measure whose density is $f$ to the measure whose
density is $g$, and is referred to as the Brenier map between the
two measures. The matrix $\nabla F(x)$ is positive-definite since it
has a positive determinant and it is the Hessian matrix of a convex function.

\medskip \textbf{Remark.}
The Knothe map, used in Section 6, is in some sense a limiting case of the Brenier map. See \cite{cgf}.

\medskip
The following lemma contains the central idea of this section.
\begin{lemma}
Let $f,g$ be log-concave probability densities in $\RR^n$. Denote $K = K(f,g)$.
Let $x \to F(x)$ be the Brenier map pushing forward  the measure whose density is  $f$ to the measure whose density is $g$. Suppose that $X$ is a random vector distributed according to the law $l(f,g)$ in $\RR^n$. Then,
\begin{equation}
Var[|X|^2] \geq \frac{1}{K(f,g)} \int_{\RR^n} f(x) Var \left [ \left |D^{-1/2}((1 - \Lambda) x + \Lambda F(x) - b(f,g)) \right |^2 \right ] dx
\end{equation}
where $D = D(f,g)$ and $\Lambda$ is a random variable distributed uniformly in $[0,1]$.
\end{lemma}

\textbf{Proof:} By a standard approximation argument we may assume that $f$ and $g$ are sufficiently smooth. Denote $D = D(f,g)$ and $L(\lambda, x) = L(f,g)(\lambda, x)$. Furthermore, define,
$$
\tilde f(x) = \sqrt{\det D} \cdot f(D^{1/2} x + b(f,g)), ~~ \tilde g(x) = \sqrt{\det D} \cdot g(D^{1/2} x + b(f,g))
$$
so that $\tilde f(x) = K(f,g) L(0, x)$ and $\tilde g(x) = K(f,g) L(1,x)$.
Denote
$$
\tilde F(x) = D^{-1/2} (F(D^{1/2}x + b(f,g)) - b(f,g)).
$$
Then $\tilde{F}$ pushes forward the measure whose density is  $\tilde{f}$   to the measure whose density is $\tilde{g}$.
Next, define
$$
M(\lambda, x) = (M_1(\lambda, x), M_2(\lambda,x)) = (\lambda, (1 - \lambda) x + \lambda \tilde F(x)).
$$
By elementary properties of the Brenier map, $M$ is a one-to-one  map
from $[0,1] \times Supp( \tilde f)$ to $Supp(L)$. Define a
density,
$$
q(\lambda, x) = \frac{\tilde f(x)^{(1 - \lambda)} \tilde g(\tilde F(x))^\lambda}{K(f,g)} = L(0,x)^{1-\lambda} L(1,\tilde{F}(x))^{\lambda}.
$$
Using the fact that $L$ is log-concave, we obtain
\begin{equation} \label{qsmall}
q(\lambda, x) \leq L(M(\lambda, x)), ~~~ \forall \lambda \in [0,1],
x \in Supp(\tilde f).
\end{equation}
A simple calculation shows that the Jacobian of $M(\lambda, x)$ is
$$
J(\lambda, x) = \det ((1 - \lambda) Id + \lambda \nabla \tilde F(x)).
$$
Recall that $\det (\nabla \tilde F(x)) = \frac{\tilde f(x)}{\tilde g(\tilde F(x))}$.
Furthermore, the matrix $\nabla \tilde F(x)$ is diagonalizable with positive eigenvalues, since it is conjugate
to the matrix $\nabla F(D^{1/2} x + b(f,g))$ which is a positive-definite matrix. By the arithmetic/geometric means inequality,
$$
J(\lambda, x) \geq \det (\nabla \tilde F(x))^{\lambda} = \left ( \frac{\tilde f(x)}{\tilde g(\tilde F(x))} \right )^\lambda.
$$
Therefore,
\begin{equation}\label{Jqsmall}
J(\lambda, x) q(\lambda, x) \geq \frac{\tilde f(x)}{K(f,g)}, ~~~\forall \lambda \in [0,1], x \in \RR^n.
\end{equation}
By changing variables using $M^{-1}$ and applying (\ref{qsmall}) and (\ref{Jqsmall}), we calculate
\begin{align*}
Var \left [ \left |X \right |^2 \right ] & =
\int_{\RR^n} \int_{[0,1]} \left ( |x|^2 - \int_{\RR^n} \int_{[0,1]} |y|^2  L(\theta, y) d \theta dy \right )^2 L(\lambda, x) d \lambda dx
\\ & \geq
\int_{\RR^n} \int_{[0,1]} \left ( |M_2(\lambda, x)|^2 - \int_{\RR^n} \int_{[0,1]}   |y|^2  L(\theta, y) d \theta dy \right )^2 J(\lambda, x) q(\lambda, x) d \lambda dx \\ & \geq
\int_{\RR^n} \frac{\tilde f(x)}{K(f,g)} \left ( \int_{[0,1]} \left ( |M_2(\lambda, x)|^2 - \int_{\RR^n} \int_{[0,1]} |y|^2  L(\theta, y) d \theta dy \right )^2  d \lambda \right ) dx \\ & \geq
\int_{\RR^n} \frac{\tilde f(x)}{K(f,g)} \left (\int_{[0,1]} \left ( |M_2(\lambda, x)|^2 - \int_{[0,1]} |M_2(\theta, x)|^2 d \theta \right )^2 d \lambda  \right ) dx \\ & =
\int_{\RR^n} \frac{\tilde f(x)}{K(f,g)}  Var \left [ \left |(1 - \Lambda) x + \Lambda \tilde F(x) \right |^2 \right ] dx.
\end{align*}
Applying the change of variables $x \to D^{-1/2} (x - b(f,g))$ completes the proof.
\qed

\medskip
By the definition of the thin-shell parameter $\sigma_n$ from \cite{EK1}, for any isotropic random vector $X$
in $\RR^n$ with a log-concave density, one has,
\begin{equation}
Var[|X|^2] \leq C n \sigma_n^2.
\end{equation}
Combining this with the above lemma yields
\begin{equation} \label{varlambda}
\int_{\RR^n} f(x) Var \left [ \left |D(f,g)^{-1/2}((1 - \Lambda) x + \Lambda F(x) - b(f,g)) \right |^2 \right ] dx \leq C K(f,g) n \sigma_n^2.
\end{equation}
For $x,y \in \RR^n$, define,
$$
v(x,y) = Var \left [|\Lambda x + (1 - \Lambda)y|^2 \right ]
$$
In view of (\ref{varlambda}), we would like to have a lower bound for $v(x,y)$ in terms of $|x|^2 - |y|^2$ and in terms of $|x-y|$.
The following lemma serves this purpose.

\begin{lemma} \label{estlambda}
There exist universal constants $C_1, C_2 > 0$, such that for all $x,y \in \RR^n$,
\begin{equation} \label{varest1}
v(x,y) = C_1 (|x|^2 - |y|^2)^2 + C_2 |x-y|^4.
\end{equation}
\end{lemma}

\textbf{Proof:} Define
$$
f(\lambda) = |\lambda x + (1 - \lambda) y|^2, ~~ g(\lambda) = \lambda |x|^2 + (1 - \lambda) |y|^2,
$$
and $h(\lambda) = f(\lambda) - g(\lambda)$. Then $h(1-\lambda) = h(\lambda)$
hence  $COV(g(\Lambda), h(\Lambda)) = 0$. Consequently,
\begin{equation} \label{eqqq1}
Var[f(\Lambda)] = Var[h(\Lambda)] + Var[g(\Lambda)].
\end{equation}
It is easy to verify that
\begin{equation} \label{eqqq2}
Var[g(\Lambda)] = (|x|^2 - |y|^2)^2 Var(\Lambda) = C_1 (|x|^2 - |y|^2)^2.
\end{equation}
Next, using the parallelogram law,
$$
h(\lambda) = -\lambda (1 - \lambda) |x-y|^2.
$$
Consequently,
\begin{equation} \label{eqqq3}
Var[h(\Lambda)] = |x-y|^4 Var \left[ \Lambda(1-\Lambda) \right] = C_2 |x-y|^4.
\end{equation}
Combining (\ref{eqqq1}), (\ref{eqqq2}) and (\ref{eqqq3}) completes the proof.
\qed

\textbf{Proof of Theorem \ref{thmsec6_half}:}
Write $b = b(f,g)$ and $D = D(f,g)$. Substituting the result of Lemma \ref{estlambda} into  (\ref{varlambda}) yields
\begin{align} \label{eqnplug}
\int_{\RR^n} f(x) & \left ( \left ( |D^{-1/2} (x-b)|^2 - |D^{-1/2} (F(x) - b)|^2 \right )^2 + |D^{-1/2} (x - F(x))|^4 \right ) dx \\
\nonumber &
\leq C K(f,g) n \sigma_n^2.
\end{align}
Let $X,Y$ be the random vectors whose densities are $f,g$ respectively. By the definition of the transportation distance,
\begin{equation} \label{eqnw2}
W_2^2(D^{-1/2}X,D^{-1/2}Y) \leq \int_{\RR^n} f(x) |D^{-1/2}(x - F(x))|^2 dx,
\end{equation}
where the transportation distance between random vectors is defined to be the distance between
the corresponding distribution measures. The fact that $f$ and $g$ have barycenters at the origin implies
$$
\EE [\langle D^{-1/2} X, D^{-1/2} d \rangle] = \EE [\langle D^{-1/2} Y, D^{-1/2} d \rangle] = 0,
$$
and consequently
\begin{align} \label{eqntrace}
\int_{\RR^n} f(x)  & \left ( |D^{-1/2} (x-d)|^2 - |D^{-1/2} (F(x) - d)|^2 \right ) dx \\ & =
 Tr(Cov(D^{-1/2}X) - Cov(D^{-1/2} Y)). \nonumber
\end{align}
The Cauchy-Schwartz inequality together with (\ref{eqnplug}), (\ref{eqnw2}) and (\ref{eqntrace}) yield,
\begin{equation} \label{finaltreq}
W_2(\tilde X, \tilde Y)^4 + \left [ Tr(Cov(\tilde X) - Cov(\tilde Y)) \right ]^2 \leq C n K(f,g) \sigma_n^2
\end{equation}
where $\tilde X = D^{-1/2} X$ and $\tilde Y = D^{-1/2} Y$.
Consequently,
$$
W_2(X,Y)^2 \leq C \sqrt{n K(f,g)} \sigma_n ||D||_{OP}
$$
where $\| D \|_{OP} = \sup_{0 \neq x} |D(x)| / |x|$ is the operator norm of $D$.
From the remark to Corollary \ref{boundedvar} we conclude
 that
$$
||D||_{OP} \leq  C K_{1/2}(f,g)^4.
$$
The function $\lambda \mapsto K_{\lambda}(f,g)$ is log-concave and it is bounded from below by one, according to the Pr\'ekopa-Leindler
inequality. Therefore,
$$ K_{1/2}(f,g) \geq \sqrt{\sup_{\lambda \in (0,1)} K_{\lambda}(f,g)} \geq \sqrt{K(f,g)}, $$
and
(\ref{thmsec6eq2}) is proven. \qed

\medskip The rest of this section aims at a better understanding of the exponents in Theorem \ref{thmsec2}. The next lemma
exploits  the second summand in our basic estimate (\ref{finaltreq}).

\begin{lemma} \label{sumitup}
Let $f,g$ be log-concave probability densities on $\RR^n$ whose barycenters are at the origin. Suppose that $f$ is isotropic. Then there exists a universal constant $c_1 > 0$ such that whenever $K_{1/2}(f,g) \leq \exp(n^{c_1})$, there exist two unit vectors $\theta_1, \theta_2 \in \RR^n$
with
\begin{equation} \label{thetaexists1}
\left  \langle Cov(g) \theta_1, \theta_1 \right \rangle  \leq 1 + C \sigma_n \sqrt{\frac{K(f,g)}{n}}
\end{equation}
and
\begin{equation} \label{thetaexists2}
\left  \langle Cov(g) \theta_2, \theta_2 \right \rangle \geq 1 - C \sigma_n \sqrt{\frac{K(f,g)}{n}}.
\end{equation}
Here,  $C > 0$ is some universal constant.
\end{lemma}

\textbf{Proof:}  We use the notation of the proof of Theorem \ref{thmsec6_half}.
In order to establish (\ref{thetaexists1}), we fix $\alpha >0$, and assume that
$$
\langle Cov(g) \theta, \theta \rangle > 1 + \alpha \sigma_n \sqrt{\frac{K(f,g)}{n}}, ~~~ \forall \theta \in \Sph,
$$
where $\Sph = \{ x \in \RR^n ; |x| = 1 \}$. Our goal is to show that necessarily $\alpha \leq C$. Noting that $Cov(\tilde X) = D^{-1}$ we have
$$
\langle Cov(\tilde X)^{-1/2} Cov(\tilde Y) Cov(\tilde X)^{-1/2} \theta, \theta \rangle - 1 > \alpha \sigma_n \sqrt{\frac{K(f,g)}{n}}, ~~~ \forall \theta \in \Sph,
$$
where $\tilde{X}$ and $\tilde{Y}$ are as in the proof of Theorem \ref{thmsec6_half}.
The last inequality implies,
$$
 \frac{Tr(Cov(\tilde Y))}{Tr(Cov(\tilde X))} - 1  > \alpha \sigma_n \sqrt{\frac{ K(f,g) }{n}}.
$$
Consequently,  in order to establish  (\ref{thetaexists1}), it suffices to show that for some universal constant $C>0$,
$$
\left | Tr(Cov(\tilde Y)) - Tr(Cov(\tilde X)) \right | \leq C Tr(Cov(\tilde X)) \sigma_n \sqrt{\frac{K(f,g) }{n}}.
$$
In view of (\ref{finaltreq}), the last inequality will be concluded if we only manage to show,
\begin{equation} \label{tracelarge}
Tr(Cov(\tilde X)) = Tr(D^{-1}) \geq \frac{n}{2}.
\end{equation}
The above fact follows from an application of Lemma \ref{lemsec5} with $\delta = 1/2$ and from the assumption that $K_{1/2}(f,g) \leq \exp(n^{c_1})$. Equation (\ref{thetaexists1}) is established, and the proof of (\ref{thetaexists2}) is analogous. The proof of the lemma is thus complete.
\qed

\medskip Next, define
\begin{equation}
\kappa = \limsup_{n \to \infty} \frac{\log \sigma_n}{\log n}, ~~~
\tau_n = \max \left \{ 1, \max_{1 \leq j \leq n} \frac{\sigma_j}{j^\kappa} \right \},
\end{equation}
so that $\sigma_n \leq \tau_n n^\kappa$. Note that the thin-shell
conjecture implies that $\kappa = 0$ and $\tau_n < C$.
We apply the estimate from the previous lemma for various marginals of our $n$-dimensional
measures, and obtain:

\begin{lemma} \label{lastlem}
Let $f,g$ be log-concave probability densities in $\RR^n$ whose barycenter is at the origin. Suppose that $f$ is isotropic. Define $R = K_{1/2}(f,g)$ and denote by $\{\lambda_i\}$ the eigenvalues of $Cov(g)$, repeated according to their multiplicity. Assume that the sequence
$\{ |\lambda_i - 1| \}$ is non-increasing. Then, one has
\begin{equation} \label{smallis}
|\lambda_i - 1| \leq C R^4, ~~ \forall 1 \leq i \leq n
\end{equation}
and
\begin{equation} \label{largeis}
|\lambda_i - 1| \leq C R \tau_n i^{\kappa - \frac{1}{2}}, ~~\forall (\log (2R))^{C_1} \leq i \leq n
\end{equation}
where $C,C_1 > 0$ are some universal constants.
\end{lemma}

\textbf{Proof:} The bound  (\ref{smallis}) follows directly from the remark to Corollary \ref{boundedvar}.
In order to establish  (\ref{largeis}), denote by $\{ e_i \}$ the orthonormal basis of eigenvectors corresponding to the eigenvalues $ \{ \lambda_i \}$. Define
$$
E_1 = sp \{e_j; ~~ 1 \leq j \leq i, \lambda_j \geq 1  \}, ~~E_2 = sp \{e_j; ~~ 1 \leq j \leq i, \lambda_j \leq 1  \}.
$$
Let $E$ be the subspace with the larger dimension among these two subspaces. Then $k = \dim E \geq i/2$.
Denote by $i_0$ the maximal $j$ for which $e_j \in E$. Then $k \leq i_0 \leq i$.
According to our assumption, $\dim(E) \geq (\log (2 R))^{C_1} / 2$, and hence we may apply Lemma \ref{sumitup}
in the subspace $E$. Denote by $f_E$ and $g_E$ the marginals of $f$ and $g$ to the subspace $E$.
Using  (\ref{thetaexists1}) and (\ref{thetaexists2})
for $f_E$ and $g_E$ we obtain
\begin{equation}
|\lambda_i - 1| \leq |\lambda_{i_0} - 1| \leq C \sigma_k \sqrt{\frac{K(f,g)}{k}} \leq C' R \tau_k i^{\kappa - \frac{1}{2}}
\leq  C' R \tau_n i^{\kappa - \frac{1}{2}}
\end{equation}
where we used the  fact that $K(f,g) \leq K_{1/2}(f,g)^2 = R^2$ as well as the Pr\'ekopa-Leindler inequality which implies  that
$K_{\lambda}(f_E, g_E) \leq K_{\lambda}(f,g)$ for any $\lambda \in (0,1)$.
\qed

\medskip The next theorem demonstrates that the exponent $\alpha_1$
in Theorem \ref{thmsec2} may be made arbitrarily close to  $1/2-\kappa$, thus complementing
the inequality (\ref{eq_353}) which goes in the opposite direction.
This provides yet another piece of evidence for the close relationship between the thin shell problem and the stability
of the Brunn-Minkowski inequality in high dimensions.

\begin{theorem} \label{thmsec6}
Let $K, T \subset \RR^n $ be convex bodies and let $R \geq 1$. Assume
that $K$ is isotropic, that the barycenter of $T$ is at the origin and that
\begin{equation}  Vol_n \left( \frac{K+T}{2} \right) \leq C R \sqrt{Vol_n(K) Vol_n(T)}. \label{eq_1235}
\end{equation}
 Then,
\begin{equation} \label{thmsec6eq1_}
\| Cov(\mu_T) - Id \|_{HS}  \leq C \left ( R^5 + \tau_n R \max(\sqrt{\log n},
n^{\kappa} ) \right ),
\end{equation}
where $Id$ is the identity matrix.
Consequently,
\begin{equation}
 \left| \frac{\int_T |x|^2 d \mu_T(x)}{\int_K |x|^2 d \mu_K(x)} - 1 \right|
\leq \frac{ \| Cov(\mu_T) - Id \|_{HS}}{\sqrt{n}} \leq  C \frac{R^5 + \tau_n R \max(\sqrt{\log n},
n^{\kappa} ) }{\sqrt{n}}. \label{eq_1238}
\end{equation}
Here, $C > 0$ is a universal constant.
\end{theorem}

\textbf{Proof:}
We may clearly assume that $Cov(\mu_T)$ is a diagonal matrix whose
diagonal is $\lambda_1,\ldots,\lambda_n$, where the sequence
$\{ |\lambda_i - 1| \}$ is non-increasing. Since our measures are log-concave, then we may use Lemma \ref{lastlem} and  calculate
\begin{align*}
\sum_{i=1}^n |\lambda_i - 1|^2 & \leq C R^8 (\log (2R))^{C_1} +  C R^2 \tau_n^2 \sum_{i=1}^n i^{2 \kappa - 1} \leq \tilde{C} R^9 +  C R^2 \tau_n^2 \left( 1 + \int_1^n s^{2 \kappa - 1} ds \right )  \\ &  \leq
C' (R^9 + \tau_n^2 R^2 \max(\log n, n^{2 \kappa} )).
\end{align*}
The bound (\ref{thmsec6eq1_}) follows.
In order to deduce (\ref{eq_1238}) from (\ref{thmsec6eq1_}), argue as in
 (\ref{eq_1234}) above.
The proof is complete. \qed

\section{Transportation in one dimension}

In this section we recall some basic definitions concerning
transportation of one-dimensional measures. For a Borel measure
$\mu$ in $\RR^n$ we write $\overline{Supp(\mu)}$ for  the set of
all points $x \in \RR^n$ such that all of the neighborhoods of $x$ have
positive $\mu$-measure. The support of $\mu$, denoted by $Supp(\mu)$, is defined in this paper
to be the interior of $\overline{Supp(\mu)}$.
Suppose that $\mu_1$ and $\mu_2$ are Borel
probability measures on the real line, with continuous densities
$\rho_1$ and $\rho_2$, respectively. We further assume that the
$Supp(\mu_2)$ is connected and that $\rho_2$ does not vanish in
its support. For $t \in \RR$ set
$$ \Phi_j(t) = \mu_j \left( (-\infty, t] \right) \quad \quad \quad \quad (j=1,2). $$
For $j=1,2$, the map $\Phi_j^{-1}$ pushes forward the uniform
measure on $[0,1]$ to $\mu_j$. The {\it monotone transportation map}
between $\mu_1$ and $\mu_2$ is the continuous, non-decreasing
function $$ F(t) = \Phi_2^{-1}(\Phi_1(t)), $$ defined for $t \in
Supp(\mu_1)$. Observe  that
$$ F_*(\mu_1) = \mu_2. $$
Furthermore, $F$ is differentiable in $Supp(\mu_1)$ and
\begin{equation}
 \rho_1(t) = F^{\prime}(t) \rho_2(F(t)) \ \ \ \ \ \ \ \ \text{for} \ \ \ t \in Supp(\mu_1).
\label{eq_1855}
\end{equation}
Additionally, it is well-known (see, e.g., Villani's book \cite{villani}) that
\begin{equation}
 W_2(\mu_1, \mu_2) = \sqrt{ \int_{\RR} |F(x) - x|^2 d \mu_1(x) }. \label{eq_1237}
 \end{equation}
A probability measure
on $\RR$ is said to be {\it even} if $\mu(A) = \mu(-A)$ for any
measurable $A \subset \RR$, where $-A = \{ -x ; x \in A \}$.

\begin{proposition} Suppose that $\mu_1$ and $\mu_2$ are even, log-concave probability
measures on $\RR$. Denote $\sigma = \sqrt{ Var(\mu_1) + Var(\mu_2)
}$.  Then,
\begin{equation} \label{compdists}
W_2(\mu_1, \mu_2) \leq C \sigma \sqrt{\int_{\RR} \min \{ (F^{\prime}(t)  - 1)^2, 1 \}  d\mu_1(t)}
\end{equation}
 where $F$ is the monotone transportation map between $\mu_1$ and $\mu_2$ and $C > 0$ is a  universal constant. \label{prop_120}
\end{proposition}

We begin the proof of Proposition \ref{prop_120} with the  following
crude lemma.

\begin{lemma} \label{wass}
Let $\mu_1$ and $\mu_2$ be probability measures on the real line.

\begin{enumerate}
\item[(i)] If $\mu_1$ and $\mu_2$ are even, then
$$
W_2(\mu_1, \mu_2)^2 \leq 2 (Var(\mu_1) + Var(\mu_2)).
$$
\item[(ii)]
If $\mu_1, \mu_2$ are supported on $[A, \infty)$ and $[B, \infty)$
respectively, and have non-increasing densities, then
$$
W_2(\mu_1, \mu_2) \leq |B-A| + 10 \sqrt{Var(\mu_1) + Var(\mu_2)}.
$$
\end{enumerate}
\end{lemma}
\textbf{Proof:} Denote by $\delta_0$ the Dirac measure at the
origin. Assume that $\mu_0$ and $\mu_1$ are even. By the triangle
inequality for the transportation metric,
$$ W_2(\mu_1, \mu_2) \leq W_2(\mu_1, \delta_0) + W_2(\delta_0, \mu_2) =
\sqrt{Var(\mu_1)} + \sqrt{Var(\mu_2)},  $$ and (i) follows. We move on
to prove (ii). Denote $e = \EE[\mu_1]$. It follows from the
fact that the density of $\mu_1$ is non-increasing that the
expectation of $\mu_1$ is larger than its median. Hence, $$ \mu_1
\left( \left[A,e \right] \right) \geq \frac{1}{2}, \quad\text{and}
\quad \mu_1 \left( \left[A,A+\frac{e-A}{2} \right] \right) \geq
\frac{1}{4}. $$ Therefore,
$$
Var(\mu_1) \geq \int_A^{A+\frac{e-A}{2}} (e-x)^2 d \mu_1(x) \geq
\frac{(e-A)^2}{16}.
$$
Let $\delta_A, \delta_B, \delta_e$ be the Dirac measures supported
on $A, B, e$ respectively. By the triangle inequality,
$$
W_2(\mu_1, \delta_A) \leq W_2(\mu_1, \delta_e) + W_2(\delta_e,
\delta_A) = \sqrt{Var(\mu_1)} + (e-A) \leq 5 \sqrt{Var(\mu_1)}.
$$
In the same manner,
$$
W_2(\mu_2, \delta_B) \leq 5 \sqrt{Var(\mu_2)}.
$$
Therefore, by using $W_2(\mu_1, \mu_2) \leq W_2(\mu_1, \delta_A) +
W_2(\delta_A, \delta_B) + W_2(\delta_B, \mu_2)$,
$$
W_2(\mu_1, \mu_2) \leq 10 \sqrt{Var(\mu_1) + Var(\mu_2)} + |B-A|.
$$
\qed

\medskip
\textbf{Proof of Proposition \ref{prop_120}:} Use (\ref{eq_1855}), the definition of $F$, and the
fact that $\Phi_1^{-1}$ pushes forward the uniform measure on
$[0,1]$ to $\mu_1$, in order to obtain
$$ \int_{\RR} \min \{ (F^{\prime}(t)  - 1)^2, 1 \}  d\mu_1(t) = \int_0^1
\min \left \{ \left(
\frac{\rho_1(\Phi_1^{-1}(t))}{\rho_2(\Phi_2^{-1}(t))} - 1 \right)^2,
1 \right \}  d t. $$ Recall that when $\mu_j$ is a log-concave
measure, the function $\rho_j(\Phi_j^{-1}(t))$ is concave on
$[0,1]$. Denote $I_j(t) = \rho_j(\Phi_j^{-1}(t))$ for $j=1,2$. Then $I_1$ and $I_2$ are
concave, non-negative functions on $[0,1]$, with the property
that $I_j(t) = I_j(1-t)$ for any $t \in [0,1]$. These two functions are
therefore continuous on $(0,1)$, increasing on $[0,1/2]$, and
decreasing on $[1/2,1]$. Let $\eps > 0$ be such that
\begin{equation}
 \eps^2 =
 \int_0^1 \min \left \{ \left( \frac{I_1(t)}{I_2(t)} - 1 \right)^2, 1 \right \}  d t.
\label{eq_1925}
\end{equation}
Suppose first that $\eps > 1/10$. In this case, from part (i) of lemma
\ref{wass},
$$
W_2(\mu_1, \mu_2)^2 \leq 2 \left( Var(\mu_1) + Var(\mu_2) \right).
$$
So whenever $\eps > 1 / 10$, the inequality (\ref{compdists}) holds
trivially for a sufficiently large universal constant $C > 0$.

\medskip From now on, we restrict attention to the case where $\eps \leq 1/10$. We divide the rest of the proof into several steps.

\medskip
\textbf{Step 1:} Let us prove that there exists a universal constant
$C>0$ such that
\begin{equation}
\int_{2 \eps^2}^{1 - 2 \eps^2} \left( \frac{I_1(t)}{I_2(t)} - 1
\right)^2 d t \leq C \eps^2. \label{eq_1325}
\end{equation}
To that end, we will show that
\begin{equation}
 I_1(t) \leq 4 I_2(t) \ \ \ \ \ \ \ \ \ \ \text{for all} \ \ t \in
[2 \eps^2, 1 - 2 \eps^2]. \label{eq_1145+}
\end{equation}
Once we prove (\ref{eq_1145+}), the desired bound (\ref{eq_1325})
 follows from (\ref{eq_1925}). We thus focus on the proof of
 (\ref{eq_1145+}). Suppose that $t_1 \in (0,1/2]$ satisfies $I_1(t_1)
> 4 I_2(t_1)$. We will
 show that in this case
\begin{equation}
 t_1 \leq 2 \eps^2.
\label{eq_1918}
\end{equation}
If $I_1(t) > 2 I_2(t)$ for all $t \in (0, t_1)$, then $t_1 \leq
\eps^2$ according to (\ref{eq_1925}). Thus (\ref{eq_1918}) holds
true in this case. Otherwise, there exists $0 < t < t_1$ with
$I_1(t) \leq 2 I_2(t)$. Let $t_0$ be the supremum over all such $t$.
Since $I_1$ and $I_2$ are continuous and non-decreasing on
$(0,t_1]$, then
$$ I_1(t_0) = 2 I_2(t_0) \leq 2 I_2(t_1) < I_1(t_1) / 2. $$
Since $I_1$ is concave, non-decreasing and non-negative on
$[0,t_1]$, then necessarily $t_0 < t_1/2$. We conclude that $I_1(t)
> 2 I_2(t)$ for any $t \in [t_1 / 2, t_1]$. From (\ref{eq_1925}) it
follows that $t_1 \leq 2 \eps^2$. Therefore (\ref{eq_1918}) is
proven in all cases. By symmetry, we conclude (\ref{eq_1145+}), and
the proof of (\ref{eq_1325}) is complete.

\medskip
\textbf{Step 2:} For any $0 \leq T \leq \Phi_1^{-1}(1 - 2 \eps^2)$
we have
$$ \int_{-T}^{T} (F^{\prime}(t)  - 1)^2  d\mu_1(t)  \leq
\int_{2 \eps^2}^{1 - 2 \eps^2} \left( \frac{I_1(t)}{I_2(t)} - 1
\right)^2 d t \leq C \eps^2, $$ where the last inequality is the
content of Step 1. Denote $\nu = \mu_1|_{[-T,T]}$, an even
log-concave probability measure. According to Lemma \ref{lem_1213},
we have $Var(\nu) \leq Var(\mu_1) \leq \sigma$. Note that the
function $F(t) - t$ is odd, hence its $\nu$-average its zero. Using
the Poincar\'e-type inequality in Lemma \ref{spectral_gap}, we see
that for any $0 \leq T \leq \Phi_1^{-1}(1 - 2 \eps^2)$,
\begin{equation}
 \int_{-T}^{T} (F(t)  - t)^2 d \mu_1(t) \leq 12 Var(\nu)
\int_{-T}^{T} (F^{\prime}(t)  - 1)^2 d\mu_1(t) \leq \tilde{C}
\sigma^2 \eps^2. \label{eq_1408}
\end{equation}

 \medskip
\textbf{Step 3:} Let $T_1 = \Phi_1^{-1}(1 - 3 \eps^2)$ and let $T_2 =
\Phi_1^{-1}(1 - 2 \eps^2)$. We use (\ref{eq_1408}) and conclude that
there exists $T_1 \leq T \leq T_2$ with
\begin{equation}  |F(T) - T|^2 \leq \tilde{C} \sigma^2 \eps^2 \, \left/ \,  \mu_1 \left( [T_1, T_2] \right) \right.= \tilde{C} \sigma^2.
\label{eq_1207} \end{equation} Denote $\nu_1 = \mu_1|_{[T, \infty)}$
and $\nu_2 = \mu_2|_{[F(T), \infty)}$. These are log-concave
probability densities with $Var(\nu_1) + Var(\nu_2) \leq \sigma^2$.
Note that we have, owing to (\ref{eq_1408}), \begin{align*}
W_2(\mu_1, \mu_2)^2 & = \int_{-T}^T (F(t)-t)^2 d \mu_1(t) + 2
\int_T^\infty (F(t)-t)^2 d \mu_1(t) \\ & \leq \tilde C \sigma^2
\eps^2 + 2 \mu_1([T,\infty)) W_2(\nu _1, \nu_2)^2.
\end{align*}
In order to prove the lemma it remains to show that $W_2(\nu _1,
\nu_2)^2 \leq C \sigma^2.$ But in view of (\ref{eq_1207}), the
latter is a direct consequence of part (ii) in lemma \ref{wass}:
Since $T, F(T) > 0$, then the log-concave densities of $\nu_1$ and
$\nu_2$ are non-increasing. This completes the proof. \qed

\medskip Let $f, g: \RR \rightarrow [0, \infty)$ be  log-concave functions with finite, positive integrals.
Denote by $\mu_f, \mu_g$ the probability measures on $\RR$ whose densities are proportional
to $f$ and $g$, respectively. Let $F$ be the monotone transportation map between $\mu_f$ and $\mu_g$.
Then $S(x) = ( F(x) + x ) / 2$ is a strictly-increasing, continuous map in $Supp(\mu_1)$. Define
\begin{equation}
h \left( S(x) \right) = \sqrt{ f(x) g(F(x)) } \quad \quad \quad \quad (x \in Supp(\mu_f)).
\label{eq_1116}
\end{equation}
We set $h(x) = 0$ for any $x$ which is not in the image
of $Supp(\mu_1)$ under $S$. Then $h$ is a well-defined, non-negative, measurable function on $\RR$.
Observe that for any $x \in \RR$,
$$ h(x) \leq \sup_{y \in \RR} \sqrt{ f(x - y) g(x + y) }. $$
We thus view the function $h$ as a refined variant of the supremum-convolution of $f$ and $g$. The following
proposition is a stability estimate for the Pr\'ekopa-Leindler inequality in one dimension. It may be viewed
as the transportation-metric version of the $L^1$-stability estimates from
 Ball and B\"or\"oczky \cite{BB}.

\begin{proposition} Suppose that $f$ and $g$ are even, log-concave functions on $\RR$ with finite, positive integrals.
Denote by $\mu_f, \mu_g$ the probability measures on $\RR$ whose densities are proportional to $f,g$
respectively. Set $\sigma = \sqrt{ Var(\mu_f) + Var(\mu_g)
}$.  Then,
\begin{equation} \label{eq_1225}
W_2^2(\mu_f, \mu_g) \leq C \sigma^2 \left( \frac{\int_{\RR} h}{\sqrt{\int_{\RR} f \int_{\RR} g}} - 1 \right)
\end{equation}
 where the function $h$ is defined via  (\ref{eq_1116}) and $C > 0$ is a  universal constant. \label{prop_1202}
\end{proposition}

\textbf{Proof:} Multiplying the functions $f$ and $g$ by positive
constants, if necessary, we may assume that $\int f = \int g = 1$.
Indeed, neither the left-hand side nor the right-hand side of
(\ref{compdists}) is changed under such normalization. Let $F$ be
the monotone transportation map between $\mu_f$ and $\mu_g$ and as
before, $S(x) = (F(x) + x) / 2$ for $x \in Supp(\mu_f)$. Applying the
change of variables $y = S(x)$ we see that
$$ \int_{\RR} h(y) dy = \int_{Supp(\mu_f)} h(S(x)) S^{\prime}(x) dx = \int_{Supp(\mu_f)} \sqrt{f(x) g(F(x))} \frac{F^{\prime}(x) + 1}{2} dx. $$
According to (\ref{eq_1855}), we have $F^{\prime}(x) g(F(x)) = f(x)$
for any $x$ in the support of $\mu_f$. Since $g$ is log-concave,
it does not vanish in $Supp(\mu_g)$, and hence $F^{\prime}(x) \neq 0$ for
any $x \in Supp(\mu_f)$. Therefore,
$$ \int_{\RR} h(y) dy = \int_{Supp(\mu_f)} \frac{F^{\prime}(x) + 1}{2 \sqrt{F^{\prime}(x)}} f(x) dx \geq
\int_{Supp(\mu_f)} \left(1 + c \min \left \{ \left( F^{\prime}(x) - 1 \right)^2, 1 \right \}  \right) f(x) dx, $$
where we used Lemma \ref{lem_1457}(ii) in the last passage. Since $\int f = 1$, then
$$ \int_{\RR} h(y) dy - 1 \geq c \int_{Supp(\mu_f)} \min \left \{ \left( F^{\prime}(x) - 1 \right)^2, 1 \right \}  f(x) dx.  $$
We may thus apply Proposition \ref{prop_120} and deduce that
$$ \int_{\RR} h(y) dy - 1 \geq c \int_{\RR} \min \left \{ \left( F^{\prime}(x) - 1 \right)^2, 1 \right \}  d \mu_f(x) \geq \frac{\tilde{c}}{\sigma^2} W_2(\mu_f, \mu_g)^2 $$
and the proposition is proven. \qed

\section{Unconditional Convex Bodies}

In this section we  prove Theorem \ref{uncond} together with its
close variant, Theorem \ref{uncond2} below. We say that a function
$\rho$ on $\RR^n$ is unconditional if it is invariant under
coordinate reflections, i.e., if $$ \rho(x_1,...,x_n) = \rho(\pm
x_1,...,\pm x_n) $$ for all $(x_1,...,x_n) \in \RR^n$ and for any
choice of signs.  For two functions $f, g: \RR^n \rightarrow [0,
\infty)$  we abbreviate \begin{equation}  H(f,g)(x) =
\sup_{y \in \RR^n} \sqrt{ f(x +  y) g (x - y) }. \label{eq_1758}
\end{equation}
Thus, $H(f,g) = H_{1/2}(f,g)$ as defined in (\ref{eq_1145}).
We will frequently consider $H(f,g)(x)$ when the functions $f$ and $g$ are defined only on a subset of
$\RR^n$. For the purpose of (\ref{eq_1758}) we treat such functions as zero outside
their original domain of definition.

\begin{theorem} \label{uncond2}
Let $M > 0$ and consider the cube $Q^n = [-M,M]^n \subset \RR^n$.
Suppose that $f,g: Q^n \rightarrow [0, \infty)$ are unconditional,
log-concave probability densities. Then,
\begin{equation}  W_2^2( \mu_f, \mu_g) \leq C M^2 \left[ \int_{Q^n} H(f,g)  -1 \right], \label{eq_1138}
\end{equation}
 where $C
> 0$ is a universal constant and $\mu_f, \mu_g$ are the probability measures with densities $f,g$
respectively.
\end{theorem}

The main tool in the
proof of Theorem \ref{uncond2} is the Knothe map from \cite{KN}, which we define next.
Let $M, f,g$ be as in Theorem \ref{uncond2}. Then
the support of $\mu_g$ is a convex set, and  $g$ does not
vanish in $Supp(\mu_g)$. The Knothe map between
$\mu_f$ and $\mu_g$ is the continuous function $F =
(F_1,\ldots,F_n): Supp(\mu_f)
 \rightarrow Supp(\mu_g)$ for which
\begin{enumerate}
\item[(a)] $\displaystyle F_*(\mu_f) = \mu_g$.
\item[(b)] For any $j$, the function $F_j(x_1,\ldots,x_n)$ actually depends  only
on the variables $x_1,\ldots,x_j$. We may thus speak of
$F_j(x_1,\ldots,x_j)$.
\item[(c)] For any $j$ and for any fixed $x_1,\ldots,x_{j-1}$, the function $F_j(x_1,\ldots,x_j)$ is non-decreasing in
$x_j$.
\end{enumerate}
It may be proven by induction on $n$ (see \cite{KN}) that the Knothe
map between $\mu_f$ and $\mu_g$ exists, and that in fact, the three
requirements above determine the function $F$ completely. Denoting
$\lambda_j(x) = \left.
\partial F_j(x) \right/ \partial x_j \geq 0$, it follows from
property (b) that
\begin{equation}  \prod_{j=1}^n \lambda_j(x) = J_F(x) =
\frac{f(x)}{g( F(x))} \label{eq_1450} \end{equation} for any $x \in Supp(\mu_1)$, where $J_F(x)$ is the Jacobian of the
map $F$. Below we
will also use the fact that the map $x \mapsto x + F(x)$, defined
for $x \in Supp(\mu_f)$, is one-to-one, as follows from properties
(b) and (c).   Set $$ \pi(x_1,\ldots,x_n) = (x_1,\ldots,x_{n-1}) $$ and let
$f_{n-1}, g_{n-1}$ be the densities of the probability measures
$\pi_*(\mu_f), \pi_*(\mu_g)$, respectively. Then $f_{n-1}$ and
$g_{n-1}$ are unconditional and log-concave. Write $T_n = F =
(F_1,\ldots,F_n)$ for the Knothe map between $\mu_f$ and $\mu_g$,
and set
$$ T_{n-1}(x_1,\ldots,x_{n-1}) = \left( F_1(x_1), F_2(x_1,x_2), \ldots, F_{n-1}(x_1,\ldots,x_{n-1}) \right). $$
Then $T_{n-1}$ is the Knothe map  between $\pi_*(\mu_f)$ and
$\pi_*(\mu_g)$. Observe that for fixed $(x_1,\ldots,x_{n-1}) \in
\pi(Supp(\mu_f))$, the map
$$ x_n \mapsto F_n(x_1,\ldots,x_n) $$
is the monotone transportation map between the probability densities
proportional to $$ t \mapsto f(x_1,\ldots,x_{n-1}, t) \quad
\text{and} \quad s \mapsto g(z_1,...,z_{n-1}, s), $$ for
$(z_1,\ldots,z_{n-1}) = T_{n-1}(x_1,\ldots,x_{n-1})$.   For $i=n-1,n$ we set
$$ S_i(x) = \frac{x + T_i(x)}{2}  $$
which is a one-to-one, continuous function, defined for $x \in Supp(\mu_f)$ when $i=n$ and for $x \in \pi \left( Supp(\mu_f) \right)$ when
$i=n-1$. According to (\ref{eq_1450}) and to property (b),
the Jacobian $J_{S_i}(x)$ of the map
$S_i$ satisfies
\begin{equation}  J_{S_i}(x) = \prod_{j=1}^i \left( \frac{1 + \lambda_j(x)}{2} \right) \geq
\prod_{j=1}^i \sqrt{ \lambda_j(x)} = \sqrt{ J_{T_i}(x) }. \label{eq_1457}
\end{equation}
Finally, for $i=n-1,n$ set
\begin{equation} \label{eq_1458}
 V(f_i,g_i) \left( S_i(x) \right) = \sqrt{ f_i( x ) g_i( T_i(x)) } \leq H(f_i,g_i) \left( S_i(x) \right).
\end{equation}
Since $S_i$ is one-to-one, then $V(f_i,g_i)$ is a well-defined
function on a subset of $Q^i$. We extend $V(f_i,
g_i)$ to the entire $Q^i$ by setting it to be zero outside its original domain of
definition.

\begin{lemma} Let $\vphi: Q^{n-1} \rightarrow [0, \infty)$ be a measurable function. Then,
$$ \int_{Q^{n-1}} \vphi(S_{n-1}(y)) f_{n-1}(y) dy \leq \int_{Q^{n-1}} \vphi(y) V(f_{n-1}, g_{n-1})(y) dy. $$
\label{lem_1445}
 \end{lemma}

 \textbf{Proof:} We use (\ref{eq_1450}) for the Knothe map $T_{n-1}$ to conclude that
 \begin{align*}
\int_{Q^{n-1}} \vphi(S_{n-1}(y)) f_{n-1}(y) dy & =
\int_{Supp(f_{n-1})} \vphi(S_{n-1}(y)) \sqrt{ f_{n-1}(y)
g_{n-1}(T_{n-1}(y)) }  \sqrt{J_{T_{n-1}}(y)} dy \\ & \leq
\int_{Supp(f_{n-1})} \vphi(S_{n-1}(y)) V(f_{n-1}, g_{n-1})(
S_{n-1}(y) )  J_{S_{n-1}}(y) dy
 \end{align*}
 where we used (\ref{eq_1457}) and (\ref{eq_1458}) in the last passage.
The map $S_{n-1}$ is one-to-one in the support of $f_{n-1}$. Changing variables $z = S_{n-1}(y)$
we obtain
$$ \int_{Q^{n-1}} \vphi(S_{n-1}(y)) f_{n-1}(y) dy \leq \int_{S_{n-1}(Supp(f_{n-1}))} \vphi(z) V(f_{n-1}, g_{n-1})(z) dz $$
and the  lemma is proven. \qed

\medskip
The following lemma will  serve as the induction step in the proof
of Theorem \ref{uncond2}.

\begin{lemma}
Let $M > 0, Q^n = [-M,M]^n$. Suppose that $f,g: Q^n \rightarrow \RR$
are unconditional, log-concave probability densities. Let $T_n, T_{n-1}, f_{n-1}, g_{n-1}$  be as above. Then,
\begin{align} \label{eq_1514}
& \int_{Q^n} |T_n(x) - x|^2 f(x) dx \\ & \leq \int_{Q^{n-1}}
|T_{n-1}(y) - y|^2 f_{n-1}(y) dy + C M^2 \left[ \int_{Q^n} V(f, g) -
\int_{Q^{n-1}} V(f_{n-1}, g_{n-1}) \right], \nonumber
\end{align}
 where $C
> 0$ is a universal constant (in fact, it is the same constant as in Proposition \ref{prop_1202}). \label{lem_1233}
\end{lemma}

\textbf{Proof:} In this proof we use $x = (y,t) \in \RR^{n-1} \times \RR$
as coordinates in $\RR^n$. From  the definition of $T_{n-1}$,
\begin{align*}
& \int_{Q^n} |T_n(x) - x|^2 f(x) dx \\ & = \int_{Q^{n-1}} |T_{n-1}(y) - y|^2 f_{n-1}(y) dy
+ \int_{-M}^M \int_{Q^{n-1}} |F_n(y, t) - t|^2  f(y,t) dy dt. \end{align*}
In order to prove the lemma, it therefore suffices to show that
\begin{equation}  \int_{-M}^M \int_{Q^{n-1}} |F_n(y, t) - t|^2  f(y,t) dy dt \leq C M^2 \left[ \int_{Q^n} V(f, g) - \int_{Q^{n-1}} V(f_{n-1}, g_{n-1}) \right]. \label{eq_1236} \end{equation}
Recall that $t \mapsto  F_n(y,t)$ is the monotone transportation map
between the even, log-concave probability measures supported on $[-M,M]$, whose densities are proportional to
$t \mapsto f(y,t)$ and $s \mapsto g(T_{n-1}(y), s)$.
The variance of an even measure supported on $[-M, M]$ cannot exceed $
M^2$. We may
therefore use Proposition \ref{prop_1202}, together with
(\ref{eq_1237}), to conclude that for any $y \in \pi(Supp(\mu_f))$,
\begin{equation}
\int_{-M}^M |F_n(y, t) - t|^2  \frac{f(y,t)}{f_{n-1}(y)} dt \leq C
M^2  \left[ \frac{\int_{-M}^M V(f,g)(S_{n-1}(y),t) dt}
{\sqrt{f_{n-1}(y) g_{n-1}(T_{n-1}(y))}}  - 1 \right].
\label{eq_1506} \end{equation}  In particular, the right-hand side of (\ref{eq_1506}) is
non-negative. We use the definition (\ref{eq_1458}) and integrate
with respect to $y$. This yields:
\begin{align*}
& \int_{Q^{n-1}} \int_{-M}^M |F_n(y, t) - t|^2  f(y,t) dt dy \leq C M^2 \int_{Q^{n-1}} \left[ \frac{\int_{-M}^M V(f,g)(S_{n-1}(y),t) dt} {V(f_{n-1}, g_{n-1})(S_{n-1}(y))}  - 1 \right]  f_{n-1}(y) dy \\ & \leq
C M^2 \int_{Q^{n-1}} \left[ \frac{\int_{-M}^M V(f,g)(y,t) dt} {V(f_{n-1}, g_{n-1})(y)}  - 1 \right]  V(f_{n-1}, g_{n-1})(y) dy
\end{align*}
where the last passage is legal according to Lemma \ref{lem_1445}.
The desired estimate (\ref{eq_1236}) follows, and the proof is complete. \qed

\textbf{Proof of Theorem \ref{uncond2}:} We will prove by induction on the dimension $n$ that
\begin{equation}  \int_{Q^n} |T_n(x) - x|^2 f(x) dx  \leq C M^2 \left[ \int_{Q^n} V(f,g)  -1 \right], \label{eq_1510}
\end{equation}
where $C$ is the constant from Lemma \ref{lem_1233}. The case $n=1$
follows from Proposition \ref{prop_1202} and from the fact that the
variance of an even measure supported on $[-M, M]$ cannot exceed $
M^2$. We assume that (\ref{eq_1510}) is proven for dimension $n-1$
and proceed with the proof for dimension $n$. Apply the induction
hypothesis for the unconditional, log-concave probability densities  $f_{n-1},
g_{n-1}$ and conclude that
\begin{equation}  \int_{Q^{n-1}} |T_{n-1}(y) - y|^2 f_{n-1}(y) dy  \leq C M^2 \left[ \int_{Q^{n-1}} V(f_{n-1},g_{n-1}) -1 \right]. \label{eq_1513}
\end{equation}
Combining (\ref{eq_1514}) and (\ref{eq_1513}),
\begin{align*} & \int_{Q^n} |T_n(x) - x|^2 f(x) dx  \\ & \leq C M^2  \left \{ \left[ \int_{Q^{n-1}} V(f_{n-1},g_{n-1})  -1 \right] + \left[ \int_{Q^n} V(f, g) - \int_{Q^{n-1}} V(f_{n-1}, g_{n-1}) \right]   \right \} \end{align*}
and (\ref{eq_1510}) is proven for dimension $n$, hence for all
dimensions. Using (\ref{eq_1510}) and the fact that $V(f,g) \leq
H(f,g)$, the theorem follows by the definition of transportation
 distance. \qed

\medskip The uniform measure on a convex body is a prime example for a log-concave
measure. Consequently, we may deduce Theorem \ref{uncond} from
Theorem \ref{uncond2} by using a crude ``cut with a big cube'' argument. The
logarithmic factor of Theorem \ref{uncond} may be an artifact
of this clumsy procedure.

\medskip \textbf{Proof of Theorem \ref{uncond}:} Let $0 \leq \gamma \leq 1/2$ be a parameter to be specified later on.
For $\alpha, \beta > 0$ we denote
$$ K_{\alpha} = K \cap [-\alpha \log n, \alpha \log n]^n, \quad \quad
T_{\beta} = T \cap [-\beta \log n, \beta \log n]^n. $$ According to
Corollary \ref{boundedvar}, we have $Cov(\mu_T) \leq CR^4$.
 Using Lemma \ref{lem_310} and a union bound, we deduce that
 \begin{equation}  \mu_K( K \setminus K_{\alpha} ) \leq  C n^{1 - c \alpha}, \quad \quad
\mu_T( T \setminus T_{\beta} ) \leq  C n^{1 - c \beta / R^2}. \label{eq_1630}
\end{equation}
We now select $\alpha$ and $\beta$ so that
$$ \mu_K( K \setminus K_{\alpha} ) = \mu_T( T \setminus T_{\beta} ) = \gamma. $$
According to (\ref{eq_1630}),
\begin{equation}
\alpha \leq C \left( 1 + \frac{\log (1/\gamma)}{\log n} \right), \quad \quad
\beta \leq C R^2 \left( 1 + \frac{\log (1/\gamma)}{\log n} \right). \label{eq_1631}
\end{equation}
Denote by $\mu_K^1$ the uniform probability measure on $K_{\alpha}$
and  similarly for $T$. By elementary properties of the
transportation metric $W_2$, it follows  that
$$
W_2^2(\mu_K, \mu_T) \leq \mu_K( K_{\alpha}) \cdot W_2^2(\mu_K^1, \mu_T^1) + \mu_K(K \setminus
K_{\alpha}) \cdot \left[ Diam(K) + Diam(T) \right]^2,$$ where $Diam(K) =
\sup_{x,y \in K} |x-y|$ is the diameter of $K$. It is well-known
(see \cite{MP}) that $Diam(K) \leq C n \sqrt{ \| Cov(\mu_K) \|_{OP}
}$ and therefore,
\begin{equation}
W_2^2(\mu_K, \mu_T) \leq W_2^2(\mu_K^1, \mu_T^1) + C \gamma n^2
R^4.\label{eq_1625} \end{equation} Note that $\mu_K^1$ and $\mu_T^1$
satisfy the requirements of Theorem \ref{uncond2} with $M = \max \{
\alpha,  \beta \} \cdot \log n$. Denote $f(x) = 1_{K_{\alpha}}(x) /
Vol_n(K_{\alpha}), g(x) = 1_{T_{\beta}}(x) / Vol_n(T_{\beta})$.
Then,
$$ \int_{\RR^n} H(f,g) = \frac{Vol_n([K_{\alpha} + T_{\beta}]/2)}{\sqrt{ Vol_n(K_{\beta}) Vol_n(T_{\beta})}} \leq \frac{R}{1 - \gamma}
\leq R( 1 + 2 \gamma) = 1 + (R - 1) + 2 R \gamma. $$
From Theorem \ref{uncond2} and (\ref{eq_1625}) we conclude that
\begin{align} \nonumber
W_2^2(\mu_K, \mu_T) & \leq C  \log^2 n \cdot [\alpha^2 + \beta^2] \cdot \left \{ (R - 1) + 2 R \gamma \right \}  + C \gamma n^2 R^4 \\
& \leq  C  \log^2 n \cdot \left[R^4 \left( 1 + \frac{\log
(1/\gamma)}{\log n} \right)^2 \right] \cdot \left \{ (R - 1) + 2 R
\gamma \right \}  + C \gamma n^2 R^4. \label{eq_1640}
\end{align}
All that remains is to select $\gamma$. In the case where $R \leq
n^{2}$, we choose $$ \gamma  = (R-1)^5 \log^2 n / (10 n^4 R^4) \leq
1/2 $$ and deduce the desired bound (\ref{eq_305}) from
(\ref{eq_1640}). In the case where $R \geq n^{2}$, we select
$\gamma = 1/2$ and still deduce (\ref{eq_305}). The theorem is thus
proven for all cases. \qed

\medskip
Next, we explain why Theorem \ref{thmsec2} provides a non-trivial
estimate for the thin-shell parameter, and why Theorem \ref{uncond}
provides yet another proof for the thin-shell estimate from
\cite{uncond}, up to logarithmic factors. Observe that when $K \subset \RR^n$ is a convex
body and $T \subset K$, then
$$ Vol_n \left( \frac{T +
K}{2} \right) \leq Vol_n(K) = R \sqrt{Vol_n(K) Vol_n(T)} $$ for $R =
\sqrt{Vol_n(K) / Vol_n(T)}$.
As before, we write $B_2^n = \{ x \in
\RR^n ; |x| \leq 1 \}$ for the Euclidean unit ball, centered at the
origin in $\RR^n$.

\begin{proposition} Let $A > 0$ and let $K \subset \RR^n$ be an isotropic convex
body. For $s > 0$ denote $K_s = K \cap (s B_2^n)$. Assume that
\begin{equation}
 \left| \, \frac{ \int_{K_s} |x|^2 d \mu_{K_s}(x)}{\int_K |x|^2 d \mu_K(x)} \ - \
 1 \,
\right|  \leq A \label{eq_334}
\end{equation}
for any $s > 0$ with $Vol_n(K_s) / Vol_n(K) \in [1/8, 7/8]$. Then,
\begin{equation}  \int_K \left( \frac{|x|^2}{n} - 1 \right)^2 d \mu_K(x) \leq C A^2
\label{eq_439} \end{equation} where $C > 0$ is a universal constant.
\label{thin_shell}
\end{proposition}

\textbf{Proof:} Standard bounds on the distribution of polynomials
on high-dimensional convex sets (see Bourgain \cite{bourgain} or
Nazarov, Sodin and Volberg \cite{nsv}) reduce the desired inequality
(\ref{eq_439}) to the estimate \begin{equation} \mu_K \left(  \left
\{  x \in K ; \left| \frac{|x|^2}{n} - 1 \right| \geq 20 A \right \}
\right) \leq \frac{1}{2}. \label{eq_440}
\end{equation}
In order to prove (\ref{eq_440}), select $a > 0$ such that
$Vol_n(K_a) = Vol_n(K) / 4$. From (\ref{eq_334}),
$$ \max_{x \in K_a} \frac{|x|^2}{n} \geq \int_{K_a} \frac{|x|^2}{n} d \mu_{K_a}(x) \geq 1 - A, $$
or equivalently,
\begin{equation}  \mu_K \left(  \left
\{  x \in K ; \frac{|x|^2}{n}   \leq 1-A \right \} \right) \leq
\frac{1}{4}. \label{eq_450}
\end{equation}
For the upper bound, let $s < t$ be such that $Vol_n(K_s) = 3
Vol_n(K) / 4$ and $Vol_n(K_t) = 7 Vol_n(K) / 8$. Then, from
(\ref{eq_334}),
\begin{align*} 1 + A & \geq \int_{K_t} \frac{|x|^2}{n} d \mu_{K_t}(x) \geq \frac{6}{7} \int_{K_s} \frac{|x|^2}{n} d
\mu_{K_s}(x) + \frac{1}{7} \max_{x \in K_s} \frac{|x|^2}{n} \\ &
\geq \frac{6}{7} (1 - A) + \frac{1}{7} \max_{x \in K_s}
\frac{|x|^2}{n}.
\end{align*}
Hence, $\max_{x \in K_s} \frac{|x|^2}{n} \leq 1 + 13 A$, or
equivalently,
\begin{equation}  \mu_K \left(  \left
\{  x \in K ; \frac{|x|^2}{n}   \geq 1 + 13 A \right \} \right) \leq
\frac{1}{4}. \label{eq_500}
\end{equation}
It is now clear that (\ref{eq_440}) follows from (\ref{eq_450}) and
(\ref{eq_500}). \qed

 {

\end{document}